\documentclass[11pt,a4paper]{amsart}
\usepackage[utf8]{inputenc}
\usepackage[T1]{fontenc}
\usepackage{lmodern}
\usepackage{amsmath,amstext,amssymb,amsopn,amsthm,enumitem,mathrsfs,mathtools,dsfont}
\usepackage[most]{tcolorbox}
\usepackage[margin=1in]{geometry}
\usepackage{cite}
\mathtoolsset{showonlyrefs}
\usepackage[colorlinks,citecolor=blue,urlcolor=blue,bookmarks=true]{hyperref}
\usepackage{tikz}
\usepackage{pgfplots}
\usepgfplotslibrary{groupplots}
\pgfplotsset{compat=1.18}
\usepackage{subcaption}
\usepackage[normalem]{ulem}

\usepackage{enumitem}
\usepackage[dvipsnames]{xcolor}
\usepackage[colorinlistoftodos,prependcaption,textsize=tiny]{todonotes}
\newtheorem{theorem}{Theorem}[section]
\newtheorem{lemma}[theorem]{Lemma}
\newtheorem{proposition}[theorem]{Proposition}
\newtheorem{corollary}[theorem]{Corollary}

\theoremstyle{definition}

\newtheorem{example}[theorem]{Example}
\newtheorem{remark}[theorem]{Remark}

\newcommand{\C}{\mathbb{C}}

\newcommand{\E}{\mathbb{E}}

\newcommand{\N}{\mathbb{N}}

\renewcommand{\P}{\mathbb{P}}
\newcommand{\R}{\mathbb{R}}
\newcommand{\Rd}{{\mathbb{R}^d}}

\newcommand{\Z}{\mathbb{Z}}

\newcommand{\PP}{\mathbf{P}}
\newcommand{\PPa}[1]{\PP^{(#1)}}
\newcommand{\EE}{\mathbf{E}}
\newcommand{\EEa}[1]{\EE^{(#1)}}

\newcommand{\calC}{\mathcal{C}}

\newcommand{\calE}{\mathcal{E}}
\newcommand{\calF}{\mathcal{F}}

\newcommand{\calH}{\mathcal{H}}

\newcommand{\calP}{\mathcal{P}}

\newcommand{\calV}{\mathcal{V}}

\newcommand{\tkappa}{{\widetilde{\kappa}}}
\newcommand{\hkappa}{{\widehat{\kappa}}}

\newcommand{\tmu}{\widetilde{\mu}}

\newcommand{\hphi}{\widehat{\phi}}

\newcommand{\phie}{\phi_\epsilon}
\newcommand{\hphie}{\hphi_\epsilon}
\newcommand{\psie}{\psi_\epsilon}
\newcommand{\phin}{\phi_n}
\newcommand{\hphin}{\hphi_n}
\newcommand{\psin}{\psi_n}
\newcommand{\xin}{\xi^{(n)}}
\newcommand{\Hn}{H^{(n)}}
\newcommand{\hH}{\widehat{H}}
\newcommand{\hHn}{\hH^{(n)}}
\newcommand{\kappan}{\kappa_n}
\newcommand{\hkappan}{\hkappa_n}

\newcommand{\hmu}{\widehat{\mu}}

\newcommand{\hW}{\widehat{W}}
\newcommand{\We}{W_\epsilon}
\newcommand{\hWe}{\hW_\epsilon}

\newcommand{\ce}{c_\epsilon}
\newcommand{\hce}{\hc_\epsilon}

\newcommand{\kappae}{\kappa_\epsilon}
\newcommand{\hkappae}{\hkappa_\epsilon}

\newcommand{\He}{H^{(\epsilon)}}
\newcommand{\xie}{\xi^{(\epsilon)}}
\newcommand{\Xe}{X^{(\epsilon)}}

\newcommand{\bkappa}{\kappa_{\beta}}
\newcommand{\bhkappa}{\hkappa_{\beta}}

\newcommand{\hc}{\widehat{c}}

\newcommand{\ga}{\mathfrak{a}}
\newcommand{\hga}{\widehat{\ga}}

\newcommand{\oga}{\overline{\ga}}
\newcommand{\hoga}{\widehat{\oga}}

\newcommand{\gu}{\mathfrak{u}}
\newcommand{\hgu}{\widehat{\gu}}

\newcommand{\rA}{{\rm A}}
\newcommand{\rAC}{{\rm AC}}

\renewcommand{\ge}{g_{\e}}
\newcommand{\gen}{\ge^{(n)}}

\newcommand{\hxi}{\widehat{\xi}}
\newcommand{\hU}{\widehat{U}}
\newcommand{\uPP}{\PP^\uparrow}
\newcommand{\uEE}{\EE^\uparrow}
\newcommand{\uPPn}{\PP^{(n),\uparrow}}
\newcommand{\uEEn}{\EE^{(n),\uparrow}}

\newcommand{\Linv}{L^{-1}}
\newcommand{\Linvn}{\big(L^{(n)}\big)^{-1}}

\newcommand{\cxi}{\check{\xi}}
\newcommand{\cA}{\check{A}}
\newcommand{\cpsi}{\check{\psi}}

\newcommand{\Rp}{\R_{+}}

\newcommand{\e}{{\rm e}}

\newcommand{\oxi}{\overline{\xi}}
\newcommand{\txi}{\widetilde{\xi}}
\newcommand{\oxin}{\oxi^{(n)}}
\newcommand{\txin}{\txi^{(n)}}

\newcommand{\tomega}{\widetilde{\omega}}
\newcommand{\tA}{\widetilde{A}}
\newcommand{\tD}{\widetilde{D}}

\newcommand{\tf}{\widetilde{f}}

\newcommand{\ed}{\stackrel{\rm d}{=}}

\newcommand{\Vo}[1]{V_{#1}}
\newcommand{\Vx}[2]{V_{#1}^{(#2)}}

\DeclareMathOperator{\supp}{supp}

\renewcommand{\Re}{\operatorname{Re}}

\newcommand{\ind}{\mathds{1}}

\newcommand{\rI}{\mathrm{I}}

\renewcommand{\subset}{\subseteq}

\newcommand{\vphi}{\varphi}
\renewcommand{\epsilon}{\varepsilon}

\renewcommand{\leq}{\leqslant}
\renewcommand{\geq}{\geqslant}

\newcommand{\ud}{\,{\rm d}}
\newcommand{\od}{{\rm d}}

\newcommand{\Csc}{\calC_c^2}
\newcommand{\Cs}{\calC^2}

\newcommand{\genxi}{\mathscr{L}}

\numberwithin{equation}{section}
\allowdisplaybreaks

\makeatletter
\providecommand\@dotsep{5}
\def\listtodoname{List of Todos}
\def\listoftodos{\@starttoc{tdo}\listtodoname}
\makeatother

\begin{document}

\title{ Moments of positive self-similar Markov processes}

\author[L. Chaumont]{Lo\"\i c Chaumont}
\address{Lo\"\i c Chaumont, Univ Angers, CNRS, LAREMA, SFR MATHSTIC, F-49000 Angers, France}
\email{loic.chaumont@univ-angers.fr}

\author[{\L{}}. Le\.{z}aj]{\L{}ukasz Le\.{z}aj}

\address{\L{}ukasz Le\.{z}aj, Wroc\l{a}w University of Science and Technology, Department of Pure and Applied Mathematics, wyb. Wyspia\'{n}skiego 27, 50-370 Wroc\l{}aw, Poland}
\email{lukasz.lezaj@pwr.edu.pl}

\date{\today}
\author[Z. Palmowski]{Zbigniew Palmowski}
\address{Zbigniew Palmowski, Wroc\l{a}w University of Science and Technology, Department of Pure and Applied Mathematics, wyb. Wyspia\'{n}skiego 27, 50-370 Wroc\l{}aw, Poland}
\email{zbigniew.palmowski@pwr.edu.pl}
\thanks{\textit{MSC2020 subject classification.} Primary 60G18; secondary 60G51.Lo\"{\i}c Chaumont is grateful for
financial support from the ANR project “Rawabranch” number ANR-23-CE40-0008.
Zbigniew Palmowski was partially supported by the National Science Centre (Poland): grant 2023/51/B/ST1/01270.
The three authors confirm that this work was carried out entirely without the use of artificial intelligence.}


\begin{abstract}
Let $X$ be a positive self-similar Markov process absorbed at the origin and let $\xi$ be the underlying L\'evy process in the Lamperti representation of $X$. The main contribution of this work is to prove the equivalence between the existence of the power moment of any order $\beta>0$ for $X$ and the existence of the exponential moment of order $\beta$ for $\xi$. 
Furthermore, we provide exact expressions for integer moments of $X$. Then given a sequence $(X^{(n)})$ of positive self-similar Markov processes issued from 0 and associated with $(\xi^{(n)})$ in the Lamperti representation, we show that for $(X^{(n)})$ to converge weakly toward $X$ in the Skohorod's topology, it is necessary and sufficient that $(\xi^{(n)})$ converges weakly toward $\xi$ and that the means of the upward ladder heights processes of $\xin$ converge to the mean of the uppward ladder height process of
$\xi$.  Some intermediate results, necessary for our proofs, consist in establishing analytical properties of Laplace exponents of the ladder height processes and identify them after exponential change of measure.
\end{abstract}

\maketitle

 %
 %

\section{Introduction}

Self-similar Markov processes occupy a central position in probability theory, often arising as the scaling limits of discrete spatial processes and playing crucial roles in branching and fragmentation mechanisms; see, e.g. Alili, Chaumont, Graczyk and \.{Z}ak \cite{17-LA-LC-PG-TZ-ejp} for review and references.

Let $(\Omega ,\mathcal{F},(\mathcal{F}_{t})_{t\geq 0},(\P_x)_{x\in(0,\infty)})$ be a filtered probability
space satisfying the usual conditions of right continuity and completion.
A standard strong Markov process $X = (X_t)_{t \geq 0}$ on $(\Omega ,\mathcal{F},(\mathcal{F}_{t})_{t\geq 0},(\P_x)_{x\in(0,\infty)})$ taking values in $\mathbb{R}_+$ is said to be a positive self-similar Markov process (pssMp) with index $\alpha > 0$ if it satisfies the scaling property: for any $c > 0$ and $x >0$, the law of the process under the measure $\mathbb{P}_{cx}$ satisfies:
\begin{equation}\label{4156}
\left( c^{-1} X_{c^\alpha t} \right)_{t \geq 0} \text{ under } \mathbb{P}_{cx} \ed (X_t)_{t \geq 0} \text{ under } \mathbb{P}_x.
\end{equation}

A foundational analytical tool for characterising these stochastic dynamics is the Lamperti representation, which provides a one-to-one correspondence between pssMp's
absorbed at 0 and time-changed L\'evy processes, see e.g.
\cite{Lamperti72, MR3160562, 17-LA-LC-PG-TZ-ejp, CC06-aopcaballero2006}. Specifically, for $t < \zeta $ (where $\zeta $ is the absorption time at 0),
\begin{equation}\label{lamperti}
X_t = x\exp \big(\xi_{\tau_{tx^{-\alpha}}} \big),
\end{equation}
where
\[
A_s = \int_0^s e^{\alpha \xi_u}\ud u, \qquad \tau_t = \inf \{s \geq 0 : A_s>t \},
\]
and $\xi$ is a possibly killed (or non conservative) L\'evy process on the space $(\Omega ,\mathcal{F})$ under some probability measure $\mathbf{P}$ (we denote by $\mathbf{E}$ the corresponding expectation). We emphasize that the transformation \eqref{lamperti} is invertible and that the paths of $\xi$ can be recovered from those of $(X_t,\,t<\zeta)$.

The primary aim of this work is to thoroughly investigate the relationship between the exponential moments of $\xi$ and the power moments of $X$. For $x,\lambda>0$, we set
\[
\Vx{\beta}{x}(\lambda)=\EE\left[\int_0^{\zeta }e^{-\lambda s}xe^{\beta \xi_{\tau_{sx^{-\alpha}}}} \ud s\right]=
\E_x\left[\int_0^{\zeta }e^{-\lambda s}X_s^\beta\ud s\right],
\]
and $\Vo{\beta}(\lambda) = V_\beta^{(1)}(\lambda)$. Note that $\Vx{\beta}{x}(\lambda)$
is the $\lambda$-potential operator of $X$ applied to the function $f(x)=x^\beta$. Then the scaling property yields:
\begin{equation}\label{eq:Vq_scaling}
\Vx{\beta}{x}(\lambda)=x^{\beta+1}\Vo{\beta}(\lambda x).
\end{equation}
In the case where the process is not absorbed at 0, a.s., conditions for the existence of a weak limit of $(\P_x)$, when $x$ tends to 0 were obtained in Caballero and Chaumont \cite{CC06-aopcaballero2006} and Chaumont et al. \cite{2012-LC-AEK-JCP-VR-aop}. When it exists, this limit is denoted by $\P_0$. 
Then the scaling property \eqref{4156} holds with $x=0$ and we set 
\[V_\beta^{(0)}(\lambda)=\E_0\left[\int_0^{\zeta }e^{-\lambda s}X_s^\beta\ud s\right].\]
Note that if $X$ is a pssMp with index $\alpha>0$ whose associated L\'evy process in the Lamperti representation \eqref{lamperti} is $\xi$, then the process $X^\alpha$ is a pssMp with index $1$ and associated L\'evy process $\alpha \xi$. This shows that our results do not depend on the value of the index $\alpha$. We shall therefore always assume $\alpha=1$. For any Borel function
$f$ we write $\EE[f(\xi_t)]$, when it exists, for $\EE[f(\xi_t),t<\e_q]=e^{-qt}\EE[f(\check{\xi}_t)]$,
where $\e_q$ is the lifetime of $\xi$ with rate $q$ and $\check{\xi}$ is the unkilled version of $\xi$.
Our first main result is as follows.
\begin{theorem}\label{thm:moments}
Assume $(X,\P_x)$ is a positive self-similar Markov process  absorbed at $0$ when $x>0$ and let $\xi$ be the (possibly killed) L\'{e}vy process associated with $X$ through the Lamperti representation
\eqref{lamperti}. Let $\beta>0$. Then the following assertions are equivalent:
\begin{enumerate}[label={{\rm (\arabic*)}}]
		\item\label{thm_moments_item_ssmp} $\mathbb{E}_x X_t^\beta<\infty$ for some (hence, all) $x>0$ and $t>0$.
		\item\label{thm_moments_item_Vx} $\Vx{\beta}{x}(\lambda)<\infty$ for some (hence, all) $\lambda>0$ and $x>0$.
		\item\label{thm_moments_item_levy} $\EE e^{\beta\xi_t}<\infty$ for some (hence, all) $t > 0$.
\end{enumerate}
Suppose furthermore that $\limsup_{t \to \infty} \xi_t=\infty$, $\PP$-a.s.~and $\xi$ is not a compound Poisson process.  
Then \ref{thm_moments_item_ssmp} --- \ref{thm_moments_item_levy}  are further equivalent to each of the following assertions,
\begin{enumerate}[label={{\rm (\arabic*)}}]
\addtocounter{enumi}{3}
		\item\label{thm_moments_item_entrance} 
        $\mathbb{P}_0$ exists and $\mathbb{E}_0 X_t^\beta<\infty$ for some (hence, all) $t>0$.
	\item\label{thm_moments_item_V0} $\mathbb{P}_0$ exists and $\Vx{\beta}{0}(\lambda)<\infty$ for some (hence, all) $\lambda>0$.
\end{enumerate}
In any case, if $\beta= n$ for some $n \in \N$, the process $\xi$ is not killed and the process $X$ is not absorbed at the origin, $\EE e^{(n+1)\xi_1}<\infty$, then  for all $t>0$ and $x>0$,
\begin{equation}\label{1882}
\mathbb{E}_x X_t^{n}=x^{n}+\sum_{k=1}^n\frac{\psi(n)\ldots\psi(n-k+1)}{k!}x^{n-k}t^k,
\end{equation}
where $\psi$ is L\'{e}vy-Khintchine exponent of $\xi$.
Relation \eqref{1882} is also true for $x=0$ whenever $\mathbb{P}_0$ exists.
\end{theorem}
Heuristically speaking, Theorem \ref{thm:moments} states that the random time change does not affect the finiteness of the power and  exponential moments of the pssMp $X$ and the L\'{e}vy process $\xi$, respectively.
We emphasise that different implications in the above statement require completely different methods in the proof. For example, the implication $\ref{thm_moments_item_ssmp} \Rightarrow \ref{thm_moments_item_levy}$ boils essentially down to a careful analysis of the infinitesimal generator of the underlying L\'{e}vy process $\xi_t$, while the direction $\ref{thm_moments_item_levy} \Rightarrow \ref{thm_moments_item_Vx}$ relies on a recursive formula for the $\lambda$-potential operator of $X_t$ applied to the function $f(x)=x^\beta$. As such, this relation (and identity \eqref{1882}) appeared in Bertoin and Yor \cite[Proposition 1]{02-JB-MY-toulouse} under stronger assumptions
and, in our setting, is reflected in Lemma \ref{lem:Vq}.
It states that for $\xi$ being non-killed, if $ e^{\psi(\beta+1)}=\EE e^{(\beta+1)\xi_1}<\infty$, then
$\Vx{\beta-1}{x}(\lambda)<\infty$ and $\Vx{\beta}{x}(\lambda)<\infty$. Moreover, 
	\begin{equation*}
	\lambda \Vx{\beta}{x}(\lambda)=\psi(\beta)\Vx{\beta-1}{x}(\lambda)+x^{\beta}.
	\end{equation*}
Inverting the appropriate Laplace transforms and applying the scaling property gives \eqref{1882}. We also note that a less explicit identity than \eqref{1882} holds also in the case where $\xi$ is killed or drifts to $-\infty$. We refer to Remark \ref{1882killed} for more details.

Let us mention the work of Rivero \cite{10.3150/07-BEJ6082}, that has strong connections with Theorem \ref{thm:moments}. Indeed Corollary 1 of the latter work asserts that when $X$ admits a recurrent extension $\tilde{X}$ that leaves 0 continuously (which implies that $\xi$ is either non conservative or drifts towards $-\infty$), for $\beta\in(0,1)$, ${\bf E}\xi_1^+e^{\beta\xi_1}<\infty$ if and only if $\mathbb{E}_x\tilde{X}_t^\beta<\infty$, for all $t,x\geq 0$. This confirms the intuition that the existence of 
a moment of order $\beta$ for the recurrent extension of $X$ implies the existence of such moments for $X$, but that the converse is false.

We also stress that the moment equivalence $\ref{thm_moments_item_levy} \Leftrightarrow \ref{thm_moments_item_ssmp}$  does not depend on neither the behaviour of $\xi$ at infinity nor the presence of the killing.
Although it is easy to see from the L\'{e}vy-Khintchine representation that the killing of $\xi$ does not affect the existence of its exponential moments (and, by Theorem \ref{thm:moments}, the power moments of $X_t$), it seems much less obvious that its behaviour at infinity also does not play any role. In particular, we allow $\xi$ to be oscillating, i.e. $\limsup_{t \to \infty}\xi_t = -\liminf_{t \to \infty} \xi_t=\infty$, $\PP$-a.s. Finally, note also that the equivalence with  $\ref{thm_moments_item_Vx}$ gives more information on the behaviour of $t\mapsto \mathbb{E}_x X_t^\beta$
(resp. $x\mapsto \mathbb{E}_x X_t^\beta$) at infinity (resp. at $0$).

\begin{remark}\label{momementscorollary}
The above main theorem has a natural consequence for a $d$-dimensional ssMp $X$.
Let $\calH$ be a locally compact subspace of $\Rd \setminus \{0\}$.
    Assume $(X,\P_x)$ is a $\calH$-valued self-similar Markov process whose norm $(|X|,\P_{|x|})$ is a positive self-similar Markov process.
    Let $(\xi,\PP)$ be a L\'{e}vy process corresponding to $(|X|,\P_{|x|})$ through the Lamperti representation. Then the following are equivalent:
    \begin{enumerate}[label={{\rm (\arabic*)}}]
    	\item $\E_x |X_t|^\beta<\infty$ for some (hence, all) $x \in \calH$ and $t > 0$.
		\item $\EE e^{\beta\xi_t}<\infty$ for some (hence, all) $t > 0$.
	\end{enumerate}
In particular, this is the case for the isotropic ssMps in the sense of \cite[Definition 3.3]{17-LA-LC-PG-TZ-ejp} and for the ssMps satisfying the skew-product property.
\end{remark}

To handle the oscillating case of $\xi$ and pssMp $X$ issued from $0$, we prove another new result that is of interest in its own.
Let $\xi$ be a conservative L\'evy process issued from $0$, satisfying $\limsup_{t\to\infty}\xi_t=\infty$, $\PP$-a.s.~and whose upward ladder height process $H$
is integrable, i.e.~${\EE}H_1<\infty$. We also assume that $\xi$ is not a compound Poisson process.
Let us now consider a sequence of L\'evy processes $(\xin)_{n\geq 0}$
satisfying the same assumptions as $\xi$ with their respective upward ladder height processes $\Hn$.
For any $n\in \mathbb{N}$, we define the pssMp $X^{(n)}$ by the Lamperti representation \eqref{lamperti} using $\xin$.
According to Chaumont at al. \cite{2012-LC-AEK-JCP-VR-aop}, we may further define pssMp's $X^{(n)}$ and $X$ starting from the origin as weak limits as $x \to 0^+$ of $X^{(n)}$ and $X$ starting from $x > 0$, respectively. We prove the following result.
\begin{theorem}\label{thm:ssMps_0_conv}
Consider a family of L\'{e}vy processes $\xin$, $n \in \N$ together with a L\'{e}vy process $\xi$, all of them conservative, satisfying $\limsup_{t\to\infty}\xi_t=\limsup_{t\to\infty}\xin_t=\infty$, $\PP$-a.s.~and whose upward ladder height processes $H$
and $\Hn$ are integrable, i.e.~${\EE}H_1<\infty$ and ${\EE}\Hn_1<\infty$. Then the sequence of  pssMp's $\Big(X^{(n)}\Big)$ issued from $0$ and corresponding to $\xin$ via the Lamperti representation \eqref{lamperti}
converges weakly in the Skorokhod topology  as $n$ tends to $\infty$ towards the pssMp $X$ issued from $0$ that corresponds to $\xi$ via the Lamperti representation if and only if the following conditions are satisfied:
\begin{enumerate}[label={\bf A\arabic*},itemindent=0.5cm, labelsep=0.5cm,itemsep=0.3cm, topsep=0.3cm]
\item\label{A:weak_conv} 
The sequence $(\xin)$ converges weakly towards $\xi$.
\item\label{A:EH_conv} 
$\lim\limits_{n\to\infty}\EE \Hn_1=\EE H_1$.
\end{enumerate}
\end{theorem}

The proof of Theorem \ref{thm:ssMps_0_conv} uses a representation analogous to that of Lamperti in \eqref{lamperti} for pssMp's issued from $0$, obtained in Bertoin and Savov \cite{2011-JB-MS-bulllms}.
The complexity of this construction makes the proof rather long and technical.
We also discuss conditions of Theorem \ref{thm:ssMps_0_conv} producing interesting counterexamples.

Finally let us mention that we also investigate how the Girsanov-type transformation affects the bivariate Laplace exponents of the ladder height processes. More precisely, we first analyse the analytic extension of the Wiener-Hopf factors associated with the L\'{e}vy process $\xi$ and then establish their representation after the exponential change of measure. These results are interesting in and of themselves, but they are also crucial tools in verifying $\ref{thm_moments_item_levy} \Rightarrow \ref{thm_moments_item_entrance}$ of Theorem \ref{thm:moments} if $\limsup_{t \to \infty} \xi_t=\infty$ $\PP$-a.s. and $\xi$ not being a compound Poisson process. In fact, in this case, using Theorem \ref{thm:ssMps_0_conv} we can shift the analysis to the transformed process $\xie$ drifting to $\infty$, for which the power moments of the corresponding pssMp issued from the origin can be accessed using the exponential functional of $\xie$, see e.g. Bertoin and Yor \cite{02-JB-MY-potanal}. Then we apply the recent results of Minchev and Savov \cite{2026-MM-MS-prob_surv} and here the identification of Wiener-Hopf factors is essential.

To our best knowledge, there are no results concerning the equivalence of power and exponential moments, as stated in Theorem \ref{thm:moments}. Bertoin and Yor \cite{02-JB-MY-toulouse} considered entire moments of positive self-similar Markov processes, whose L\'{e}vy processes in the Lamperti representation have finite exponential moments of arbitrary positive order. \cite[Proposition 1(i)]{02-JB-MY-toulouse} gives the same expression for the integer moments of $X$ in a much more restrictive setting. We are also not aware of any results concerning weak convergence of pssMp's issued from the origin, as stated in Theorem \ref{thm:ssMps_0_conv}. The analytical properties of Wiener-Hopf factors are partially analysed in Vigon's PhD thesis \cite{Vigon02_PhD}. The convergence of the Wiener-Hopf factors under the assumption of weak convergence of underlying L\'{e}vy processes was investigated in Ivanovs \cite{2018-Ivanovs-ann_appl_probab} although here we provide a detailed argument.

\subsection{Organisation of this paper}
The paper is organised as follows. In the next section, we introduce basic notations. In Section \ref{sec:h_ex}
we analyse the analytical extension of the Wiener-Hopf factors and identify the Wiener-Hopf factors after exponential change of measure. In doing so, we produce certain results that are of interest in their own, see e.g.~Theorem \ref{thm:WH_changed}.
In Section \ref{sec:weakconv} we prove the weak convergence of pssMp's started from the origin, i.e.~Theorem \ref{thm:ssMps_0_conv} and discuss the assumptions which appear in its formulation.
Finally, Section \ref{sec:moments_proof} is devoted to the proof of the main result, Theorem \ref{thm:moments}.

\section{Basic notions}

Before going into the proof of the main result and presenting all other key results in the following sections, we introduce basic notation and facts. Let $\C_{(a,b)}=\{z \in \C \colon a<\Re z < b\}$, $\C_a = \{z \in \C \colon \Re z=a\}$ and $\C_{[a,b)} = \C_{(a,b)} \cup \C_a$. The sets $\C_{(a,b]}$ and $\C_{[a,b]}$ are defined accordingly. We then let ${\rm A}_{(a,b)}$ be the set of analytic functions on $\C_{(a,b)}$ and let $\rA_{[a,b)}$ be the subset of ${\rm A}_{(a,b)}$ of analytic functions which are continuous up to $\C_a$. The sets $\rA_{(a,b]}$ and $\rA_{[a,b]}$ are defined accordingly. Similarly, ${\rm M}_{(a,b)}$ is the set of meromorphic functions on $\C_{(a,b)}$ with appropriate extensions to ${\rm M}_{(a,b]}$ and ${\rm M}_{[a,b)}$. We also write $\C^2_{(a,b)}$ for $\C_{(a,b)} \times \C_{(a,b)}$ and $\rA^2_{(a,b)}$ for the class of bivariate analytic functions on $\C^2_{(a,b)}$.

We recall that under $\PP$, the process $\xi$ is a L\'evy process appearing in the Lamperti transformation of $X$ given in \eqref{lamperti}. 
We shall also denote by  $\PP_x$, $x\in\mathbb{R}$, the probability measure under which $\xi_0=x$, almost surely.
We allow $\xi$ to be killed at the independent exponential time $\e_q$ with parameter $q \geq 0$ by setting $\xi_t=-\infty$ for $t \geq \e_q$, with the convention $\e_0=\infty$ a.s., which corresponds to the conservative (i.e. non-killed) case. It is then well known that its law is characterised by 
the characteristic exponent, 
to wit, we have
\begin{equation}\label{eq:L-K}
\EE e^{z\xi_t} = e^{t\psi(z)}, \quad z \in i\R,
\end{equation}
and $\psi \colon i\R \mapsto \C$ is given by the L\'{e}vy-Khintchine representation
\[
\psi(z) = \frac{\sigma^2}{2}z^2+\gamma z + \int_{\R} \big( e^{zx}-1- zx\ind_{\{|x|<1\}}\big)\,\nu(\od x)-q,
\]
where $\sigma^2 \geq 0$, $\gamma \in \R$ and $\nu$ is the L\'{e}vy measure satisfying
\[
\nu(\{0\})=0 \qquad \text{and} \qquad \int_\R \big(1 \wedge |x|^2\big)\,\nu(\od x)<\infty.
\]
We then consider the analytic continuation of $\psi$. To this end, let
\begin{equation*}
\rAC = \big\{z \in \C \colon \EE\big|e^{z\xi_1}\big|<\infty\big\}.
\end{equation*}
By \cite[Theorem 25.17]{Sato99}, we necessarily have $\rAC = \{z \in \C \colon \Re z \in C\}$, where $C$ is a convex set on the real line that contains the origin, so that $\psi$ and both sides of \eqref{eq:L-K} are continuous on $\rAC$ and analytic in the interior of $\rAC$, and
\begin{equation}\label{eq:Laplace_ex}
\EE e^{z\xi_t} = e^{t\psi(z)}, \quad z \in \rAC.
\end{equation}
We also note that $\psi(z) \in \R$ for $z \in C$.

Throughout this paper, $\beta$ will denote some positive constant. Assume that $\EE e^{\beta\xi_1}<\infty$, or equivalently (see \cite[Theorem 25.17]{Sato99}) that $\int_{|x| \geq 1}e^{\beta x}\,\nu(\od x)<\infty$. Then it follows that $\C_{[0,\beta]} \subseteq \rAC$ and the quantity $\psi(\beta)$ is well defined.
Hence, for every $\alpha \in (0,\beta]$
the process $e^{\alpha \xi_t-\psi(\alpha) t}$ is a $\PP$-martingale and we can consider the exponential change of measure, also known as the Esscher transform,
\begin{equation}\label{eq:cPP}
\frac{\od \PPa{\alpha}}{\od \PP}\Bigg|_{\calF_t} = e^{\alpha \xi_t-\psi(\alpha)t}.
\end{equation}
One verifies that under $\PPa{\alpha}$, the process $\xi$ is again a L\'{e}vy process with the L\'{e}vy-Khintchine exponent
\begin{equation*}
\psi_\alpha(z) = \psi(z+\alpha)-\psi(\alpha), \quad z \in \C_{[-\alpha,\beta-\alpha]},
\end{equation*}
and its characteristic triplet $(\sigma_\alpha,\gamma_\alpha,\nu_\alpha)$ is given by $\sigma_\alpha = \sigma$, $\gamma_\alpha = \gamma+\sigma^2\alpha + \int_{|x|<1}x(e^{\alpha x}-1)\,\nu(\od x)$ and $\nu_\alpha(\od x) = e^{\alpha x}\,\nu(\od x)$. Let us stress that $\xi$ under $\PPa{\alpha}$ is conservative, i.e. $\psi_\alpha(0)=0$, irrespectively of the presence of killing under $\PP$. We also have
\begin{equation}\label{eq:36}
\PPa{\alpha}(\xi_t \in \od x) = e^{-\psi(\alpha)t}e^{\alpha x}\PP(\xi_t \in \od x).
\end{equation}
Let us comment on some properties of the Esscher transform:
\begin{remark}\label{rem:Esscher_transform}
    Assume $\EE e^{\beta \xi_1}<\infty$ for some $\beta>0$.
    \begin{enumerate}
        \item If $\psi(0)=-q=0$, then $\psi_\alpha(z) \to \psi(z)$ for $z \in i\R$ as $\alpha \to 0^+$. Hence, in this case $\xi^{(\alpha)}$ converges to $\xi$ in the sense of finite-dimensional distributions, and therefore, also in the Skorokhod space (see e.g. \cite[Corollary VII.3.6]{JS87}).
        \item It is standard to verify that $\EE \xi_1 \in [-\infty,\infty)$. Under the additional assumption $\limsup_{t \to \infty} \xi_t=\infty$, $\PP$-a.s. one obtains using the law of large numbers (see \cite[Theorems 36.5 and 36.6]{Sato99}) that in fact $\EE \xi_1 \in [0,\infty)$.
        \item Differentiating \eqref{eq:Laplace_ex} for real $z \in (0,\beta)$ and taking the limit as $z \to 0^+$ gives $(\psi)'(0^+)=\EE \xi_1$.
        Differentiating \eqref{eq:Laplace_ex} twice shows the strict convexity of $[0,\beta]\ni z \mapsto \psi(z)$.
If $(\psi)'(0+) <0$ we let $\alpha_0 >0$ to be the solution of the equation
\begin{equation*}
(\psi)'(\alpha_0)=0,
\end{equation*}
if it exist for $\alpha_0<\beta$.
If $(\psi)'(0+) \geq 0$, then we choose $\alpha_0=0$.
Thus
$\EEa{\alpha} \xi_1 = \psi_\alpha'(0+)=(\psi)'(\alpha)>0$
for
        \begin{equation}\label{choiceofalpha}
        \alpha > \alpha_0,
        \end{equation}
        so that, again by the law of large numbers, $\lim_{t \to \infty} \xi_t=\infty$ $\PPa{\alpha}$-a.s. for $\alpha$ satisfying \eqref{choiceofalpha}.
    \end{enumerate}
\end{remark}

\section{Analytical extension of the Wiener-Hopf factorisation}\label{sec:h_ex}

Our goal in this section is to express the bivariate Laplace exponents of the ladder heights of $\xi$ under $\PPa{\beta}$ defined formally in \eqref{eq:cPP} in terms of those under $\PP$. Although our main result, Theorem \ref{thm:WH_changed}, seems to be intuitively predictable, we were unable to find a proper reference, so we provide a detailed proof. Let us start with a brief introduction of all the relevant objects.

Throughout this section we assume that $\xi$ is conservative, that is, $\psi(0)=0$. We then consider its ascending and descending ladder processes $(L_t^{-1},H_t)$, where $L_t$ is the local time at $0$ of $\xi$ reflected at its past supremum and $H_t = \xi_{L_t^{-1}}$ with $L^{-1}$ being the right-continuous inverse of $L$. It is then well known that $(L_t^{-1},H_t)$ is a bivariate subordinator whose Laplace exponent $\kappa$, i.e.
\[
\EE e^{-qL_t^{-1}-zH_t} = e^{-t\kappa(q,z)}, \quad q,z\geq 0,
\]
is given by,
\begin{equation}\label{eq:biv_Bernstein_f}
\kappa(q,z) = d_1q+d_2z+\int_{(0,\infty)^2} \big( 1-e^{-qt-zx}\big)\,\mu(\od t,\od x),
\end{equation}
with $d_1,d_2 \geq 0$ and
\begin{equation*}
\int_{(0,\infty)^2} \big(1 \wedge \sqrt{t^2+x^2} \big) \,\mu(\od t,\od x)<\infty.
\end{equation*}
The analogous object for the dual process $\hxi=-\xi$ will be denoted by $\hkappa$. Since $\kappa$ and $\hkappa$ are Laplace exponents of bivariate subordinators, they are naturally extended to functions which are continuous on $\C_{[0,\infty)} \times \C_{[0,\infty)}$ and analytic on each variable, that is, $\kappa,\hkappa \in \rA^2_{[0,\infty)}$. Recall that the ladder processes of the L\'evy processes $\xi$ and $\hxi$ are defined up
to a constant that depends on the normalisation of the local time in their running supremum. As is generally the practice, we 
fix here this normalisation such that $\kappa(1,0)=1$ (resp. 
$\hkappa(1,0)=1$). Then
from \cite[Corollary VI.10, p. 165]{Bertoin96}, the Laplace exponents $\kappa$ and $\hkappa$ can be identified in terms of the law of $\xi$ and $\hxi$ by the Fristedt formula \cite[Chapter VI, Corollary 10]{Bertoin96}:
\begin{align}\label{eq:Fristedt}
\begin{aligned}
\kappa(q,z) &= \exp \bigg( \int_0^\infty \int_{[0,\infty)} \big(e^{-t}-e^{-qt-zx}\big) t^{-1} \,\PP(\xi_t \in \od x) \ud t \bigg), \quad q \geq 0,\,z \in \C_{[0,\infty)},\\
\hkappa(q,z) &= \exp \bigg( \int_0^\infty \int_{[0,\infty)} \big(e^{-t}-e^{-qt-zx}\big) t^{-1} \,\PP(\hxi_t \in \od x) \ud t \bigg), \quad q \geq 0,\,z \in \C_{[0,\infty)}.\end{aligned}
\end{align}
The celebrated Wiener-Hopf factorisation (see, e.g. \cite{Bertoin96} or \cite{Kyprianou14}) then reads as
\begin{equation}\label{eq:WH}
\psi(z)-q = -h (q)\kappa(q,-z)\hkappa(q,z), \quad q \geq 0, \,z \in i\R,
\end{equation}
where
\begin{equation}\label{hq}
h(q) = \exp \bigg( -\int_0^\infty \big( e^{-t}-e^{-qt} \big)t^{-1}\PP(\xi_t=0)\ud t \bigg).
\end{equation}
We also observe that $h \colon [0,\infty) \mapsto \Rp$ is a decreasing function, which is not identically equal to 1 if and only if $\xi$ is a compound Poisson process.
First, let us establish the analytical properties of Laplace exponents $\kappa$ and $\hkappa$. A similar result appeared in Vigon's dissertation \cite{Vigon02_PhD}:
\begin{proposition}\label{prop:kappa_extension}
    Assume $\EE e^{\beta \xi_1}<\infty$. Then, for every $q \geq 0$ we have $\kappa(q,\,\cdot\,) \in \rA_{[-\beta,\infty)}$.
\end{proposition}
\begin{proof}
    Assume first that $q>0$ and define
    \[
    \vphi_q(z) = \frac{q-\psi(-z)}{h(q)\hkappa(q,-z)}, \quad z \in \C_{[-\beta,0]}.
    \]
    By assumption, $\psi \in \rA_{[0,\beta]}$, while $\hkappa(q,\,\cdot\,) \in \rA_{[0,\infty)}$. Moreover, since $q>0$, the latter function has no zeros on $\C_{[0,\infty)}$. Thus, $\vphi_q \in \rA_{[-\beta,0]}$ and, by the Wiener-Hopf factorisation on the real line \eqref{eq:WH}, it agrees with $\kappa(q,\,\cdot\,)$ on the imaginary line $i\R$. Thus, defining the extension
    \begin{align}\label{eq:kappa_ex}
    \tkappa(q,z) = \begin{cases}
        \kappa(q,z), & \Re z \geq 0,\\
        \frac{q-\psi(-z)}{h(q)\hkappa(q,-z)}, & -\beta \leq \Re z <0,
    \end{cases}
    \end{align}
    and using Morera's theorem one can prove that $\tkappa(q,\,\cdot\,) \in \rA_{[-\beta,\infty)}$.

    For $q=0$ the same approach can be applied to prove that $\tkappa(0,\,\cdot\,)$ is analytic on $\C_{(-\beta,\infty)} \setminus \{z \in i\R : \hkappa(0,z)=0\}$ and continuous up to the boundary $-\beta+i\R$, which is insufficient for our purposes. Instead, note that, for every $q>0$ and $z \in \C_{[0,\infty)}$,
    \begin{align}\label{eq:fq}\begin{aligned}
    f_q(z) \coloneqq \EE e^{-qL_1^{-1}-zH_1} &= \int_{(0,\infty)^2}e^{-qt-zx}\,\PP(L_1^{-1} \in \od t,H_1 \in \od x) \\ &=\int_{(0,\infty)} e^{-zx} \ud F_q(x),\end{aligned}
    \end{align}
    where
    \[
    F_q(x) = \int_{(0,\infty) \times [0,x]} e^{-qt}\, \PP(L_1^{-1} \in \od t,H_1 \in \od y).
    \]
    Thus, we see that $f_q(z)$ is a Laplace-Steljtjes transform of a non-decreasing function $F_q$. Recall from the first part of the proof that $f_q \in \rA_{[-\beta,\infty)}$. Therefore, by \cite[Chapter II, Theorem 5b]{Widder46}, the integral defining exponential moments in \eqref{eq:fq} is absolutely convergent for $\Re z > -\beta$, and since the LHS of \eqref{eq:fq} is continuous up to the boundary $\Re z = -\beta$, by the monotone convergence one concludes that the integral in fact converges absolutely for $\Re z \geq -\beta$. Consequently,

    \[
    \EE e^{-qL_1^{-1}+\beta H_1} = \exp (-\tkappa(q,-\beta)) = \exp\Bigg(-\frac{\psi(\beta)-q}{h(q)\hkappa(q,\beta)}\Bigg)<\infty.
    \]
    By the monotone convergence and continuity of $\psi$ and $\hkappa$, taking a limit as $q \to 0^+$ yields that
    $
    \EE e^{\beta H_1} 
    <\infty.
    $
    Thus, by \cite[Theorem 25.17]{Sato99}, exponential moments $\EE e^{-zH_1}$ exist for all $z \in \C_{[-\beta,\infty]}$ and $\tkappa(0,\,\cdot\,) \in \rA_{[-\beta,\infty)}$, proving the claim for $q=0$. One also observes that \eqref{eq:kappa_ex} holds in this case as well.
\end{proof}
In what follows, we will identify $\kappa(q,\,\cdot\,)$ with its analytic continuation $\tkappa$ from the above proof.
\begin{remark}\label{rem:H_levy_measure}
    Note that the proof of the above proposition shows in particular that $\EE e^{\beta H_1}<\infty$. Thus, if we denote its Laplace exponent by $\phi$, so that
    \begin{equation*}
    \EE e^{-\lambda H_1}=e^{-\phi(\lambda)} \quad \text{with} \quad \phi(\lambda) = b \lambda + \int_{(0,\infty)} \big(1-e^{-\lambda x}\big) \,\mu_H(\od x), \quad \lambda \geq 0,
    \end{equation*}
    then \cite[Theorem 25.17]{Sato99} entails that
    \[
    \int_{(1,\infty)} e^{\beta x}\,\mu_H(\od x) < \infty,
    \]
    and the integral defining $\phi$ is absolutely convergent for $\Re \lambda \geq -\beta$. Comparing with \eqref{eq:biv_Bernstein_f} yields $\phi(z)=\kappa(0,z)$, so that $b=d_2$ and $\mu_H(\od x)=\int_{(0,\infty)} \,\mu(\od t,\od x)$.
\end{remark}

The next proposition establishes the analytic continuation of the exponent in the Fristedt representation \eqref{eq:Fristedt}. For $\Re z \geq 0$ the result is well known and requires no additional assumptions (see e.g. \cite{Bertoin96}), but since this extension is essential for us, we choose to provide a detailed argument. 

Suppose that $\EE ^{\beta \xi_1}=e^{\psi(\beta)}<\infty$, let $\psi_* = \inf\{\psi(\lambda):\lambda \in [0,\beta]\}$ and define $\alpha_*$ by $\psi(\alpha_*)=\psi_*$. Then we define $\psi^{-1} \colon [\psi_*,\infty) \mapsto \R$ by
\[
\psi^{-1}(q) = \alpha_* \vee \begin{cases}
    \sup\{\lambda \in [0,\beta]:\psi(\lambda)=q\},& \text{if }q \in \psi([\alpha_*,\beta]),\\
    \beta,& \text{otherwise}.
\end{cases}
\]
We note that $\psi^{-1}(q) \in [\alpha_*,\beta]$ for all $q \in (\psi_*,\infty)$. The number $\psi^{-1}(q)$ will serve in a forthcoming proposition as the boundary for the analytic continuation of the RHS of \eqref{eq:Fristedt}.
Denote by $f(q,z)$ the complex integral that appears in the Fristedt representation \eqref{eq:Fristedt} of $\kappa$ , i.e.,
\[
f(q,z) = \int_0^\infty \int_{[0,\infty)} \left( e^{-t}-e^{-qt-zx}\right)t^{-1}\,\PP(\xi_t \in \od x)\ud t, \quad q \geq 0,\,\Re z \geq 0.
\]
\begin{proposition}\label{prop:Fristedt}
Assume $\EE e^{\beta \xi_1}<\infty$. Then for any fixed $q \in (\psi_*,\infty)$ we have $f(q,\,\cdot\,) \in \rA_{(-\psi^{-1}(q),\infty)}$. If, moreover, $q>\psi(\beta)$, then $f(q,\,\cdot\,)\in\rA_{[-\beta,\infty)}$.
\end{proposition}
\begin{proof}
    Fix $q \in (\psi_*,\infty)$ and let $z \in E$ for some compact subset $E$ of $\C$, so that there is $M_1>0$ and $M_2 \in \R$ such that $|z|\leq M_1$ and $\Re z \geq M_2$. We then bound the integrand of $f(q,\,\cdot\,)$ as follows:
    \begin{enumerate}[label={\rm \arabic*.}]
        \item For $t \geq 1$ and $x \geq 0$, we write
        \[
        \big|e^{-t}-e^{-qt-zx}\big| \leq e^{-t}+e^{-qt}e^{-\Re z \cdot x} \leq e^{-t} +e^{-qt}\max \left\{ e^{-M_2x},1\right\}.
        \]
        \item For $t,x \in [0,1]$, we write
        \begin{align*}
        \big|e^{-t}-e^{-qt-zx}\big| &\leq \big|e^{-t}-1\big| + \big|e^{-qt-zx} -1\big| \leq \max\left\{1,e^{q+M_1}\right\}\big((q+1)t+M_1|x|\big).
        \end{align*}
        \item For $t \in [0,1]$ and $x \geq 1$, we write
        \begin{align*}
            \big|e^{-t}-e^{-qt-zx}\big| &\leq \big|e^{-t}-1\big| + \big|e^{-qt-zx} -1\big| \leq t + 1+e^{-\Re z \cdot x} \leq t+ 2\max\left\{e^{-M_2 \cdot x}, 1\right\}.
        \end{align*}
    \end{enumerate}
    Recall that $\EE [|\varsigma_t|\wedge 1] = O(\sqrt{t})$ as $t \to 0^+$ for any L\'{e}vy process $\varsigma$, which is already enough to show the absolute convergence of the integrand of $\kappa$ in the case $M_2\geq 0$, cf. \cite[bounds in the proof of Lemma 7]{Bertoin96}. Suppose now that $-\beta\leq-\psi^{-1}(q) < M_2<0$, so that in particular $ E\subset \C_{[-\beta,\infty)}$. Using the exponential change of measure $\PPa{-M_2}$ (which is valid due to the fact that $M_2 \geq -\beta$ and $\EE e^{\beta \xi_1}<\infty$),
    \[
    \frac{\od \PPa{-M_2}}{\od \PP} \Bigg|_{\calF_t} = e^{-\psi(-M_2)\cdot t} e^{{-M_2} \cdot \xi_t},
    \]
    cf. \eqref{eq:cPP}, and the fact that $\xi$ remains a L\'{e}vy process under $\PPa{-M_2}$, we get that, for $t \in [0,1]$,
    \[
    \EE \big[e^{-M_2 \cdot \xi_t} \ind_{\{\xi_t \geq 1\}}\big] = \EEa{-M_2} \big[e^{\psi(-M_2)t}\ind_{\{\xi_t \geq 1\}} \big] \leq c\PPa{-M_2}(\xi_t \geq 1) \leq c\EEa{-M_2}\big[|\xi_t| \wedge 1\big],
    \]
    with $c=\max \big\{e^{\psi(-M_2)},1\big\}$. Hence $\EE \big[e^{-M_2 \cdot \xi_t} \ind_{\{\xi_t \geq 1\}}\big] = O(\sqrt{t})$ as $t \to 0^+$ by the same argument as above. Since $\EE e^{\beta \xi_1}<\infty$, using the bounds 1., 2.~and 3., one concludes that the integral defining $f(q,\,\cdot\,)$ is absolutely convergent for $z \in \C_{(-\psi^{-1}(q),\infty)}$. Since the above estimates are uniform in $z \in E$, the analyticity of $f(q,\,\cdot\,)$ follows by applying Morera's theorem. If $q>\psi(\beta)$, then the same arguments yield the continuity up to the boundary $\Re z=-\beta$.
\end{proof}
\begin{remark}\,
    \begin{enumerate}
        \item If $\xi$ drifts to $-\infty$, then for every $q \in (\psi_*,\infty)$ one has $f(q,\,\cdot\,) \in \rA_{[-\alpha_0,\infty)}$.
        \item From the above proof it follows that $f(q,\,\cdot\,)$, and consequently also $\kappa(q,\,\cdot\,)$, are zero-free on $\C_{(-\psi^{-1}(q),\infty)}$.  We stress, however, that Proposition \ref{prop:kappa_extension} gives $\kappa(q,\,\cdot\,) \in \rA_{[-\beta,\infty)}$, while the analytic extension of the Fristedt formula is valid only for $z \in \C_{(-\psi^{-1}(q),\infty)}$.
    \end{enumerate}
\end{remark}
We are now ready to prove Theorem \ref{thm:WH_changed}. To our best knowledge, the only mention of this kind of result can be found in \cite{Vigon02_PhD}, see Proposition 4.4.1 therein, but the assumed hypothesis {\bf H} translates to the condition $\psi(\beta)<0$ being equivalent to $\EE e^{\beta \xi_1} \in (0,1)$, which is too restrictive in our setting. For the sake of completeness and integrity, we note that a related result concerning the uniqueness of the spatial Wiener-Hopf factorisation can be found in a recent paper by D\"{o}ring, Savov, Trottner and Watson \cite{MR4808554} and a similar statement for the killed process by Kuznetsov \cite[Theorem 1(f)]{2011-AK-tpa}, but since they are different from our setting, we chose to provide a rigorous proof. Recall that every object pertaining to $\xi$ under $\PPa{\beta}$ is denoted by adding a subscript ,,$\beta$'', so that, for example, $\bkappa$ stands for the Laplace exponent of the ascending ladder process under $\PPa{\beta}$, and so on.

\begin{theorem}\label{thm:WH_changed}
    Assume that $\EE e^{\beta \xi_1}<\infty$. 
    There are constants $c_\beta,\hc_\beta \in (0,\infty)$ such that the following identifications hold for all $q > 0$ and all $z \in \C_{[0,\infty)}$:
    \[
    \bkappa(q,z) = c_\beta\kappa(q+\lambda,z-\beta) \qquad \text{and} \qquad \bhkappa(q,z)=\hc_\beta\hkappa(q+\lambda,z+\beta),
    \]
    where $\lambda=\psi(\beta)$. If additionally $\lambda \geq 0$, then the above identities also apply for $q=0$.
    Moreover, $c_\beta$ and $\hc_\beta$ can be calculated by the following formulae:
    \begin{equation}\label{eq:cp}
    c_\beta = \frac{\bkappa(q_0,0)}{\kappa(q_0+\lambda,-\beta)} = \exp \bigg( \int_0^\infty \int_{[0,\infty)} e^{-t}\big(e^{-\lambda t+\beta x}-1\big)t^{-1}\,\PP(\xi_t \in \od x)\od t \bigg),
    \end{equation}
    and
    \begin{equation}\label{eq:cm}
    \hc_\beta = \frac{\bhkappa(q_0,0)}{\hkappa(q_0+\lambda,\beta)} = \exp \bigg( \int_0^\infty \int_{(-\infty,0]} e^{-t}\big(e^{-\lambda t+\beta x}-1\big)t^{-1}\,\PP(\xi_t \in \od x)\od t \bigg),
    \end{equation}
    where $q_0>0$ is arbitrary. The ratios defining $c_\beta$ and $\hc_\beta$ are in fact independent of $q_0$ and one has $c_\beta \hc_\beta=h(\lambda)/h_\beta(0)$ for the function $h$ defined in \eqref{hq}. In particular, $c_\beta \hc_\beta=1$ whenever $\xi$ under $\PP$ is not a compound Poisson process.
\end{theorem}
\begin{proof}
    Assume first that $q>0$, so that $\psi^{-1}(q+\lambda)>\beta$. By Proposition \ref{prop:Fristedt}  we see that the integrals that define $\kappa(q+\lambda,-\beta)$ and $\kappa(q+\lambda,z-\beta)$ in the Fristedt formula \eqref{eq:Fristedt} are well defined for every $z \in \C_{[0,\infty)}$. Thus,
    \begin{align*}
    \frac{\kappa(q+\lambda,-\beta)}{\kappa(q+\lambda,z-\beta)} = &\exp\bigg( \int_0^\infty t^{-1}e^{-qt}\int_{[0,\infty)}(e^{-z x}-1)e^{-\lambda t}e^{\beta x}\PP(\xi_t \in \od x) \ud t \bigg) \\ = &\exp\bigg( \int_0^\infty t^{-1}e^{-qt}\int_{[0,\infty)}(e^{-z x}-1)\PPa{\beta}(\xi_t \in \od x) \ud t \bigg) \\ = &\frac{\bkappa(q,0)}{\bkappa(q,z)},
    \end{align*}
    where in the second equality we use \eqref{eq:36}.
    Thus,
    \[
    \bkappa(q,z) = c_{\beta,q} \kappa(q+\lambda,z-\beta), \quad z \in \C_{[0,\infty)},
    \]
    where
    \begin{equation*}
    c_{\beta,q} = \frac{\bkappa(q,0)}{\kappa(q+\lambda,-\beta)}.
    \end{equation*}
     We further observe with the help of \eqref{eq:36} that
    \begin{equation*}
    c_{\beta,q} = \exp \bigg( \int_0^\infty \int_{[0,\infty)} e^{-t}\big(e^{-\lambda t+\beta x}-1\big)t^{-1}\,\PP(\xi_t \in \od x)\od t \bigg),
    \end{equation*}
    so that the integral defining $c_{\beta,q}$ is in fact (absolutely convergent and) independent of $q$. Thus, writing $c_\beta$ instead of $c_{\beta,q}$ we may rewrite the equation above as
    \begin{equation}\label{eq:qgeq0}
    \bkappa(q,z) = c_\beta \kappa(q+\lambda,z-\beta), \quad q>0,\,z \in \C_{[0,\infty)},
    \end{equation}
    where $c_\beta \in (0,\infty)$.

    Now we extend the identity above to include $q=0$. Assume $\lambda \geq 0$ and let $\Re z \geq \beta$. Then both sides of \eqref{eq:qgeq0} are continuous as $q \to 0$ for every $z \in \C_{[\beta,\infty)}$. Thus,
    \begin{equation}\label{eq:qeq0}
    \bkappa(0,z) = c_\beta \kappa(\lambda,z-\beta), \quad z \in \C_{[\beta,\infty)}.
    \end{equation}
    By Proposition \ref{prop:kappa_extension}, both sides of \eqref{eq:qeq0} are in $\rA_{[0,\infty)}$, so by the identity theorem, \eqref{eq:qeq0} holds for all $z \in \C_{[0,\infty)}$. The argument for the dual counterpart is essentially the same.

   To finish the proof, we note that for any fixed $q_0 > 0$ we have from \eqref{hq}, \eqref{eq:cp}, \eqref{eq:cm}, the classical Wiener-Hopf factorisation \eqref{eq:WH} and its extension given by \eqref{eq:kappa_ex} that
    \[
    c_\beta \hc_\beta = \frac{h(q_0+\lambda)}{h_\beta(q_0)} = \exp \bigg( -\int_0^\infty e^{-t}\big(1-e^{-\lambda t}\big)t^{-1}\PP(\xi_t=0)\ud t \bigg) = \frac{h(\lambda)}{h_\beta(0)}.
    \]
    In particular, $c_\beta\hc_\beta=1$ whenever $\xi$ under $\PP$ is not a compound Poisson process.
\end{proof}

The following supplement will be useful:
\begin{corollary}\label{cor:WH_changed}
    Assume that $\EE e^{\beta \xi_1}<\infty$ and $\limsup_{t \to \infty}\xi_t=\infty$ $\PP$-a.s. Then for every $\gamma \in (0,\beta)$ there are constants $c_\gamma,\hc_\gamma \in (0,\infty)$ such that
    \begin{align*}
    \kappa_\gamma(q,z) &= c_\gamma \kappa(q+\lambda_\gamma,z-\gamma), \quad q \geq 0,\,z \in \C_{[-\beta+\gamma,\infty)},\\
    \hkappa_\gamma(q,z) &= \hc_\gamma \hkappa(q+\lambda_\gamma,z+\gamma), \quad q \geq 0,\,z \in \C_{[-\beta+\gamma,\infty)},
    \end{align*}
    where $\lambda_\gamma = \psi(\gamma)$. Here $c_\gamma$ and $\hc_\gamma$ are given by \eqref{eq:cp} and \eqref{eq:cm}, respectively, $c_\gamma \hc_\gamma=h(\lambda_\gamma)/h_\gamma(0)$ and $c_\gamma \hc_\gamma=1$ whenever $\xi$ under $\PP$ is not a compound Poisson process.
\end{corollary}
\begin{proof}
    Since $\EE e^{\gamma \xi_1}<\infty$ and $\psi(\gamma)>0$, by Theorem \ref{thm:WH_changed}, we get the claim for $q \geq 0$ and $z \in \C_{[0,\infty)}$. But from Proposition \ref{prop:kappa_extension}, both sides of the identities above are in $\rA_{[-\beta+\gamma,\infty)}$, so the corollary follows by the uniqueness of the analytic continuation.
\end{proof}
\begin{remark}
    Let us conclude with the observation that the results of this section have the natural counterpart under the assumption $\EE e^{-\beta \xi_1}<\infty$. Essentially, the necessary changes in all statements and proofs boil down to replacing all objects pertaining to the ascending ladder height process $H_t$ with its dual counterpart $\widehat{H}_t$.
\end{remark}


\section{Weak convergence of pssMps started from the origin}\label{sec:weakconv}
The aim of this section is to prove Theorem \ref{thm:ssMps_0_conv} and to discuss the assumptions which appear, perhaps surprisingly, in its formulation.

\subsection{Theorem \ref{thm:ssMps_0_conv} and its proof}

    First, let us recall the construction due to Bertoin and Savov \cite{2011-JB-MS-bulllms} of positive self-similar Markov processes issued from 0. In what follows, $\xi$ is a L\'evy process with infinite lifetime, satisfying $\limsup_{t\to\infty}\xi_t=\infty$, $\PP$-a.s.~and whose upward ladder height process $H$
is integrable, i.e.~${\EE}H_1<\infty$. We also assume that $\xi$ is not a compound Poisson process. Under these assumptions, the family of probability measures $(\P_x)$ converges weakly as $x\to 0^+$ towards a probability measure $\P_0$ under which $X$ satisfies $\P_0(X_0=0)=1$ and the transition probabilities
of $X$ under $\P_0$ are the same as under $\P_x$, for any $x>0$, see \cite{2012-LC-AEK-JCP-VR-aop}.\\

We will denote by $\hxi=-\xi$ the dual L\'evy process and by $\hU$ the renewal function of the
upward ladder hight process of $\hxi$. Note that under our assumptions, the latter is a possibly killed subordinator whereas $H$ is not killed. Let also $a\geq0$ be the drift coefficient of the subordinator $H$ and recall that $\nu$ is the L\'evy measure of $\xi$. Then let us introduce the measure $\rho$ on $[0,\infty)^2$ by
\begin{equation}\label{eq:rho_def}
\rho(\od x,\od y)=\frac1{{\EE}H_1}(\hU(x)\,\nu(x+\od y)\ud x+a\delta_0(\od x)\,\delta_0(\od y)).
\end{equation}
The measure $\rho$ is actually a probability measure that corresponds to the weak limit of the
undershoot-overshoot distribution ${\PP}(\xi_{\tau^+_z-}\in \od x,\,\xi_{\tau^+_z}\in \od y)$ when $z$ tends to $\infty$, where
\[\tau_z^+=\inf\{t:\xi_t>z\}.\]

\begin{remark}\label{8033}
Note that in \cite{2011-JB-MS-bulllms} the case where $\xi$ is a subordinator is not excluded. However, the elements of potential theory and fluctuation theory needed to construct the self-similar process issued from $0$ in this case are not clearly defined. Actually, all the results in \cite{2011-JB-MS-bulllms} are still true in this case provided that one considers $\widehat{H}$ as the zero subordinator, so that $\hU\equiv 1$. Then the upward ladder
height process of $\xi$ is the process $\xi$ itself and the measure $\rho$ becomes
\[
\rho(\od x,\od y)=\frac1{{\EE}\xi_1}(\nu(x+\od y)\ud x+a\delta_0(\od x)\,\delta_0(\od y)).
\]
\end{remark}

Let us now make a few reminders on L\'evy processes conditioned to stay positive. We refer to
\cite{2005-LC-RD-ejp} for more details. By definition, the law of any L\'evy process $(\xi,{\PP})$ conditioned to
stay positive and issued from $x>0$ is defined for all $t\geq0$ and all $\mathcal{F}_t$-measurable set
$\Lambda$, by
\[
\uPP_{x}(\Lambda)=
\frac1{\hU(x)}{\EE}_{x}\left[\hU(\xi_t)\ind_\Lambda\ind_{\{t<\tau_{0}^-\}}\right],
\]
where we set
\[\tau_{y}^-=\inf\{t:\xi_t<y\},\quad\text{for all $y\in\R$.}\]
It is plain that $\hU(x)>0$
for all $x>0$. Note that when $0$ is not regular for $(-\infty,0)$ under ${\PP}$, $\hU(0)>0$ and the
above definition makes sense for $x=0$. If $0$ is regular for $(0,\infty)$ under ${\PP}$, then it is proved
in \cite[Theorem 2]{2005-LC-RD-ejp}, that $\uPP_{x}$ converges weakly on $D$, when $x$ tends to 0,
toward a Markovian probability measure that we will also denote by $\uPP_{0}$. In any case,
a meaning can be given to the law $\uPP_{0}$ and $\uPP_{0}(\xi_0=0)=1$ if
and only if $0$ is regular for $(0,\infty)$ under ${\PP}$. If $0$ is not regular for $(0,\infty)$ under
${\PP}$, then $\uPP_{0}(\xi_0>0)=1$, see \cite{2005-LC-RD-ejp}. Moreover, the law
$\uPP_0$ can always be represented as the law of the post-minimum of the conditioned process
starting from any positive level, as proved in \cite[Theorem 1]{2005-LC-RD-ejp}. More formally, the time
$m=\inf\{t:\xi_t=\inf_{s\geq0}\xi_s\}$ when the process reaches its overall infimum is almost surely
unique under $\uPP_x$, for any $x>0$. Moreover, under $\uPP_x$, the
post-infimum process, $(\xi_{t+m}-\xi_{m}\wedge\xi_{m-},t\geq0)$ has law
$\uPP_{0}$. The latter representation of $\uPP_{0}$ is the one we will use in what follows.
Recall also that under the assumption $\limsup_{t \to \infty} \xi_t=\infty$ $\PP$-a.s., the process $\xi$ under $\uPP_x$ is conservative and
$\uPP_x(\lim_{t\to\infty}\xi_t=\infty)=1$ for all $x\geq0$.


Let $\tD\coloneq D(\R,\R)$ and $D \coloneq D([0,\infty),\R)$ be the sets of real valued c\`adl\`ag paths indexed by $\R$ and $[0,\infty)$, respectively.
Similarly, we denote $D_t\coloneq D([0,t],\R)$. Then we define the probability measure $\calP$ on $D$ as the law of a process $(\txi_t)_{t\in\R}$ such that for any $(x,y)\in[0,\infty)^2$, conditionally on $(-\txi_{0-},\txi_0)=(x,y)$, the processes $(-\txi_{-t-})_{t\geq0}$ and $(\txi_t)_{t\geq0}$ are independent with laws $\uPP_x$ and $\PP_y$, respectively. Moreover, the law of the pair $(-\txi_{0-},\txi_0)$ is $\rho(\od x,\od y)$. More formally, if we define the canonical bijection $\vphi\colon\tomega\mapsto (\omega,\omega')$ from $\tD$ to $D\times D$,  by
\begin{equation}\label{bijection}
\tomega(t)=\left\{\begin{array}{ll}
\omega'(t)&\text{ if }t\geq0,\\
-\omega(-t-)&\text{ if }t<0,
\end{array}\right.
\end{equation}
then the distribution $\calP$ of $\txi$ is characterised by
\begin{equation}\label{1258}
\vphi\circ\calP(\od\tomega)=
\int_{[0,\infty)^2}\,\rho(\od x,\od y)\,\uPP_x (\od\omega)\,\PP_y(\od\omega').
\end{equation}
Note that $\txi_t<0$ for all $t<0$ and $\lim_{t\to-\infty}\txi_t=-\infty$, a.s.
Furthermore, the process $(\txi_t-\txi_0,\,t\geq0)$ has the same law as the L\'evy process $\xi$.\\

Such a process $\txi$ with law $\calP$ on $D$ being given, we set for all $t\in\R$,
\[
\tA_t\coloneq\int_{-\infty}^te^{\txi_s}\od s \quad \text{and} \quad \sigma_t\coloneq\inf\{s:\tA_s>t\}.
\]
Our assumptions imply that for every $t\in\R$, $\tA_t$ is finite a.s.~and that
$\lim_{t\to\infty}\tA_t=\infty$, a.s. Then according to \cite{2011-JB-MS-bulllms}, the process
$X$ defined by
\begin{equation}\label{eq:pssMp_ex}
X_t=\exp({\txi_{\sigma_t}}),\quad t\geq0,
\end{equation}
is a self-similar Markov process (with scaling index 1) such that $X_0=0$ and $X_t>0$,
for all $t>0$, almost surely. Moreover, the transition probabilities of $X$ are the same as those of
the canonical process $X$ under $\P_x$, for all $x>0$. Its law is the weak limit of $\uPP_x$ as $x$ tends to 0.\\

Let us now consider a sequence of conservative L\'evy processes $(\xin)_{n\geq 0}$ under the measure ${\PP}$ and denote by $\Hn$ their respective upward ladder height processes. We assume that
none of the processes $\xin$ is a compound Poisson process and that
$\limsup_{t\rightarrow\infty}\xi_t=\infty$ and $\limsup_{t\rightarrow\infty}\xin_t=\infty$, a.s., for all $n \in \N$ (but we emphasise that these assumptions do not hold in
Lemma \ref{lem:WH_factors_conv} below).

Let us denote by $X^{(n)}$ the pssMp issued from $0$ that is defined relatively to $\xin$ in the same way as in \eqref{eq:pssMp_ex}.
Assume that $\EE \Hn_1<\infty$, $n \in \N$, and $\EE H_1<\infty$.
In this section, we prove Theorem \ref{thm:ssMps_0_conv}, that is, that
the sequence of processes $\Big(X^{(n)}\Big)$ converges weakly in the Skorokhod topology of $D$ towards the process
$X$, as $n$ tends to $\infty$, if and only if \ref{A:weak_conv} and \ref{A:EH_conv} hold. 

\begin{remark}\label{rem:representation_thm} Since the Skorokhod's space $D$ is Polish, it follows that the support of the
distribution of $\xi$ is separable. Hence, under assumption \ref{A:weak_conv} the Skorokhod representation theorem applies, 
that is, we can define a sequence $\big(\xi^{*,(n)}\big)$ and a process $\xi^*$ in some probability space such that
$\xi^*\ed\xi$ and for all $n\geq0$, $\xi^{*,(n)}\ed\xin$ and such that the sequence
$(\xi^{*,(n)})$ converges almost surely in the topology of $D$ towards $\xi^*$. We shall make use
of this result several times in the sequel.
\end{remark}

In order to prove Theorem \ref{thm:ssMps_0_conv} we need a couple of preliminary results. The first lemma is quite natural and has already been used in the mathematical literature, see for instance \cite[step (i) in Appendix C]{2018-Ivanovs-ann_appl_probab}. However, we have not found a detailed demonstration of this result and since it is essential in our arguments, it is worth to provide one here.

\begin{lemma}\label{lem:WH_factors_conv} Let $\xin$, $n\in \N$, and $\xi$ be L\'evy processes such
that \ref{A:weak_conv} holds and $\xi$ is not a compound Poisson process. Let $\Linvn$ and $\Linv$ be the upward ladder time processes
of $\xin$ and $\xi$, respectively, and denote by $\kappan$ and $\kappa$ the Laplace exponents
of the bi-variate subordinators $\Big(\Linvn,\Hn\Big)$ and $\big(\Linv,H\big)$. Then, for all $(\alpha,\beta)\neq(0,0)$,
\[
\lim_{n\to\infty}\kappan(\alpha,\beta)=\kappa(\alpha,\beta),
\]
and $\kappan(0,0)$ does not tend to $\kappa(0,0)$ if and only if $\xi$ drifts to $-\infty$ and $\limsup_{t \to \infty} \xin_t=+\infty$ a.s. for infinitely many $n$. Consequently, $\Big(\Linvn,\Hn\Big)$ converges towards $\big(\Linv,H\big)$ in the sense of finite-dimensional distributions if and only if $\limsup_{t \to \infty} \xi_t=+\infty$ a.s. or $\lim_{t \to \infty} \xin_t=-\infty$ a.s. from a certain rank.
\end{lemma}
\begin{remark}
    If $\limsup_{t \to \infty} \xi_t=+\infty$ a.s. and $\limsup_{t \to \infty} \xin_t=+\infty$ a.s. from a certain rank, then $\Big(\Linvn,\Hn\Big)$ in fact converges weakly towards $\big(\Linv,H\big)$ in the Skorokhod topology as $n$ tends to $\infty$.
\end{remark}
\begin{proof}[Proof of Lemma \ref{lem:WH_factors_conv}]
Let ${\e}$ be an exponentially distributed random variable with parameter $q>0$ that is independent of
$\xi$ and the sequence $(\xin)_{n\geq0}$ and define
$\gen=\sup\{s<\e:\xin_s=\oxin_s\}$, where
$\oxin_s=\sup_{u\le s}\xin_u$. The variables $\ge$ and
$\oxi_s$ are defined in the same way with respect to $\xi$. Then the distribution
of $(\gen,\oxin_\e)$ is related to the Laplace exponent of
$\Big(\Linvn,\Hn\Big)$, as follows:
\begin{equation}\label{eq:WH_Bertoin}
{\EE}\left[e^{-\alpha \gen-\beta\oxin_\e} \right]=
\frac{\kappan(q,0)}{\kappan(\alpha+q,\beta)},\quad \alpha,\beta\geq0,
\end{equation}
see \cite[Eq. (1), p. 163]{Bertoin96}. If, moreover, $\xin$ drifts toward $-\infty$,
then $\oxin_\infty:=\sup_{u\geq0}\xin_u<\infty$ and
$g_\infty^{(n)}:=\sup\{s:\xin_s=\oxin_\infty\}<\infty$, a.s. In this case, the
subordinator $\Big(\Linvn,\Hn\Big)$ is killed so that
$\kappan(0,0)>0$ and \eqref{eq:WH_Bertoin} can be extended to $q=0$ by dominated convergence and continuity
of $\kappan$ on $\R_+^2$, and we obtain
\begin{equation}\label{eq:WH_Bertoin_q0}
\EE \left[e^{-\alpha g_\infty^{(n)}-\beta\oxin_\infty}\right]=
\frac{\kappan(0,0)}{\kappan(\alpha,\beta)},\quad \alpha,\beta\geq0.
\end{equation} The same relation holds between the distribution of
$(\ge,\oxi_\e)$, $(g_\infty,\oxi_\infty)$ and the Laplace exponent $\kappa$ of $(\tau,H)$.

Now let $t>0$ and define the functions $\omega\mapsto\overline{\omega}_t=\sup_{s\le t}\omega_s$ and
$\omega\mapsto\gamma_t(\omega)=\sup\{s<t:\omega_s=\overline{\omega}_s\}$ on the set $D_t$ of c\`adl\`ag paths from
$[0,t]$ to $\R$. Observe that for $\omega\in D_t$, if $\gamma_t(\omega)$ is the unique time $s\in[0,t]$
such that $\omega_{s}=\overline{\omega}_t$ or $\omega_{s-}=\overline{\omega}_t$ then the mapping
$\omega \mapsto (\gamma_t(\omega),\overline{\omega}_t)$ is continuous 
in the topology of the uniform convergence. Moreover, if $\chi\colon [0,t]\rightarrow[0,t]$ is a continuous increasing bijection,
then $(\overline{\omega\circ\chi})_t=\overline{\omega}_t$ and $\gamma_t(\omega\circ\chi)=\chi^{-1}\circ\gamma_t(\omega)$.
Therefore, the mapping $\omega\mapsto (\gamma_t(\omega),\overline{\omega}_t)$ is also continuous 
in the Skorokhod's topology of $D_t$ (note that the continuity of the function $\omega \mapsto \overline{\omega}_t$ on $D_t$ is proved in \cite[Theorems 6.2(i) and 2.8]{1980_Whitt}).

Then recall that for all $t>0$, $g_t$ is a.s.~the unique time $s\in[0,t]$ such that
$\xi_s=\oxi_t$ or $\xi_{s-}=\oxi_t$, see \cite[Section VI.2]{Bertoin96}. Therefore, it follows from the result obtained in the former paragraph and the weak convergence of $(\xin)$
toward $\xi$ on each space $D_t$, $t>0$ (see also \cite[Theorem 2.8]{1980_Whitt}), that for all $t>0$, the
$\R_+^2$-valued random variables $\Big(g_t^{(n)},\oxin_t\Big)$ converge weakly toward
$\big(g_t,\oxi_t\big)$, as $n$ tends to $\infty$. Then $\e$ being independent of $(\xin)$ and $\xi$,
we obtain (using dominated convergence) that the $\R_+^2$-valued random variables
$\Big(g_{\e}^{(n)},\oxin_{\e}\Big)$ converge weakly toward $(g_{\e},\oxi_{\e})$,
as $n$ tends to $\infty$. We derive from \eqref{eq:WH_Bertoin} that
\begin{equation}\label{8331}
\lim_{n\to\infty}\frac{\kappan(q,0)}{\kappan(\alpha+q,\beta)}=
\frac{\kappa(q,0)}{\kappa(\alpha+q,\beta)}, \quad \alpha,\beta\geq0,\,q>0.
\end{equation}

Let us now recall from \eqref{eq:Fristedt} the Fristedt formula 
for the Laplace exponent of $(\tau,H)$ that reads for $(\tau^{(n)},H^{(n)})$ as, 
\begin{equation*}
\kappan(\alpha,\beta)=\exp\left(\int_0^\infty \ud t\int_{[0,\infty)}(e^{-t}-e^{-\alpha t-\beta x})t^{-1}{\PP}(\xin_t \in dx)\ud t\right),\;\;\;\alpha,\beta\geq 0.
\end{equation*}
For $\alpha>0$ and $\beta=0$, this leads to 
$\kappan(\alpha,0)=\exp\left(\int_0^\infty(e^{-t}-e^{-\alpha t})t^{-1}{\PP}(\xin_t \geq 0)\ud t\right)\,.$
Since $\PP(\xi_t\neq\xi_{t-})=0$ for all $t\geq0$, \ref{A:weak_conv} entails that $(\xin_t)$
converges weakly toward $\xi_t$. Moreover, since $\xi$ is not a compound Poisson process, for all $t>0$, the
distribution function of $\xi_t$ is continuous on $\R$, so that in particular,
\[
\lim_{n\to\infty}\PP(\xin_t\geq0)=\PP(\xi_t\geq0).
\]
Hence, it follows from the above expression of $\kappan(\alpha,0)$ and the dominated convergence that $\lim_{n\rightarrow0}\kappan(q,0)=\kappa(q,0)$ for
all $q>0$. We then derive from \eqref{8331} that for all $\alpha>0$ and $\beta\geq0$,
$\lim_{n\to\infty}\kappan(\alpha,\beta)=\kappa(\alpha,\beta)$. Thus, it remains to deal with the case $\alpha=0$.

Assume first that $\kappa(0,0)=0$ (or equivalently that $\xi$ does not drift to $-\infty$).
If none of the processes $\xin$, $n\in \N$, drifts to $-\infty$, that is $\kappan(0,0)=0$, then
our result is proved. If all processes $\xin$, $n\geq0$ drift to $-\infty$, then for all $n\geq0$,
$\oxin_\infty=\sup_{s\geq0}\xi_s^{(n)}<\infty$,
and $g_\infty^{(n)}=\sup\{s<\infty:\xi_s^{(n)}=\oxi_\infty^{(n)}\}<\infty$,
a.s. We readily derive
from $\lim_{t\to\infty}\oxi_t=\infty$, a.s.~and the weak convergence of $(\xin)$ toward
$\xi$ that $\lim_{n\to\infty}\oxin_\infty=\infty$ and $\lim_{n \to \infty} g^{(n)}_\infty=\infty$, weakly. Then from \eqref{eq:WH_Bertoin_q0} and
the fact that for all $\alpha>0$ and $\beta\geq0$,
$\lim_{n\to0}\kappan(\alpha,\beta)=\kappa(\alpha,\beta)>0$, we obtain that
$\lim_{n\to\infty}\kappan(0,0)=\kappa(0,0)=0$. By considering two subsequences, we see that
the same conclusion holds if some processes $\xin$ drift towards $-\infty$ and others do not.

Denote by $\mu_n$ and $\mu$ the distribution of the pair $\Big(\Linvn,\Hn\Big)$ and $\big(\Linv,H\big)$, respectively. Then from what we have already proved, the Laplace transforms $f_n(\alpha,\beta)$ of $\mu_n$ converge to the Laplace transform $f(\alpha,\beta)$ of $\mu$ for all $\alpha > 0$ and $\beta \geq 0$. The convergence of $\kappa_n(0,0)$ to $\kappa(0,0)$ translates to the convergence of masses $f_n(0,0)=\mu_n([0,\infty)^2)$ to $f(0,0)=\mu([0,\infty)^2)=1$. Thus, one may rescale $\mu_n$ to probability measures $\tmu$ and conclude that the Laplace transforms $\tf_n(\alpha,\beta)$ of $\tmu_n$ also converge to $f(\alpha,\beta)$ for all $\alpha>0$ and $\beta \geq 0$. It is then standard to verify that the family $\tmu_n$ is tight and converges weakly to $\mu$, so that its Laplace transforms $\tf_n(\alpha,\beta)$ converge to $f(\alpha,\beta)$ for all $\alpha,\beta \geq 0$, and consequently, so do $f_n(\alpha,\beta)$. It follows readily that $\kappa_n(\alpha,\beta)$ tends towards $\kappa(\alpha,\beta)$ for all $\alpha,\beta \geq 0$, as desired.

Now suppose $\kappa(0,0)>0$, that is, $\xi$ drifts to $-\infty$. Then we observe that $\hkappa(0,0)=0$ due to the fact that $\xi$ has infinite life-time. Thus, by the previous steps, we conclude that $\hkappa_n(\alpha,\beta)$ tends towards $\hkappa(\alpha,\beta)$ for all $\alpha,\beta \geq 0$. In particular, if we denote by $\hmu_n$ and $\hmu$ the distributions of $\hHn$ and $\hH$, respectively, then by rescaling $\hmu_n$ to probability measures $\widetilde\hmu_n$ as in the above paragraph, we see that $\widetilde\hmu_n$ converges properly to $\hmu$ in the terminology of Feller, see \cite[Chapter VIII.1]{feller1971introduction} for a definition and \cite[XIII.1, Theorem 2]{feller1971introduction} for a statement. Consequently, by the continuity theorem \cite[XV.3, Theorem 2]{feller1971introduction} and the convergence of $\hkappan(0,0)$ to $\hkappa(0,0)=1$ we conclude using the rescaling argument that $\hkappan(0,ia)$ tends towards $\hkappa(0,ia)$ for all $a \in \R$. Since $\xin$ can possibly live on a lattice, let $Z=\cup_{n=1}^\infty Z^{(n)}$, where $Z^{(n)}=\{a \in \R: \psin(ia)=0\}$. It is clear that $Z$ is at most countable. Since $\Hn$ and $\hHn$ live on the positive part of the same lattice as $\xin$, it follows from the Wiener-Hopf factorisation that $\kappan(0,ia)$ tends towards $\kappa(0,ia)$ for all $a \in \R \setminus Z$. But then one may check the proof of the lemma in\cite[XV.3]{feller1971introduction} and verify that the convergence almost everywhere (so in particular: for all $a \in \R \setminus \{0\}$ if all $\xin$ are non-lattice, so that $Z=\{0\}$) is enough to conclude \cite[XV.3, Equation (3.4)]{feller1971introduction}, and consequently, that the distribution $\mu_n$ of $\Hn$ converges to the distribution $\mu$ of $H$ in the sense of \cite[Chapter VIII.1]{feller1971introduction}. Then using again \cite[XIII.1, Theorem 2]{feller1971introduction} one gets that $\kappan(0,\beta)$ tends towards $\kappa(0,\beta)$ for all $\beta>0$.

Finally, let us consider $\beta=0$. Assume that by some rank $n_0$,
all processes $\xin$, $n\geq n_0$ drift to $-\infty$. Then $\oxi_\infty<\infty$ and
$g_\infty<\infty$, a.s. and for all $n\geq n_0$,
$\oxin_\infty<\infty$ and $g_\infty^{(n)}<\infty$, a.s.
Moreover, in a similar way as we proved above that the random variables
$(g_{\rm e}^{(n)},\oxin_{\rm e})$ converge weakly toward $(g_{\rm e},\oxi_{\rm e})$,
we can prove that the random variables $(g_\infty^{(n)},\oxin_\infty)$ converge weakly
toward $(g_\infty,\oxi_\infty)$, as $n$ tends to $\infty$
(indeed, consider the mapping $\omega\mapsto(\overline{\omega}_\infty,\gamma_\infty(\omega))$, where
$\overline{\omega}_\infty=\sup_{s\geq0}\omega_s$ and
$\gamma_\infty(\omega)=\sup\{s\geq0:\omega_s=\overline{\omega}_\infty\}$, on the set
$D^\downarrow\coloneqq\{\omega\in D:\lim_{t\to\infty}\omega_t=-\infty\}$ endowed with the induced
topology). Then the relation \eqref{eq:WH_Bertoin_q0} allows us to conclude that
$\lim_{n\to\infty}\kappan(0,0)=\kappa(0,0)$. Conversely, if $\limsup_{t \to \infty}\xin_t=+\infty$ for infinitely many $n$, then there is a subsequence $n_k$ such that $\kappa_{n_k}(0,0)=0$ for all $k \in \N$ and consequently, the convergence of $\kappan(0,0)$ to $\kappa(0,0)$ cannot hold.

\end{proof}

\begin{remark}
Let us provide an example illustrating the necessity of the additional assumption regarding the case where $\xi$ drifts to $-\infty$. Let $\xin:=n^{-1}S+\xi$, where $S$ and $\xi$ are independent
L\'evy processes. We assume, moreover, that $S$ is increasing, $\xi$ drifts to $-\infty$, $\EE S_1=\infty$, and $\EE |\xi_1|<\infty$.
Then for all $n\geq1$, $\lim_{t\to\infty}\xi_t^{(n)}=\infty$, a.s. Moreover, for all $t$,
$\lim_{n\to\infty}\xin_t=\xi_t$, a.s., so that the sequence $(\xin)$ converges
weakly on $D$ towards $\xi$. Therefore, $\kappan(0,0)=0$, for all $n\geq1$, whereas $\kappa(0,0)>0$.
However, as proved above, for all $(\alpha,\beta)\neq (0,0)$ one has $\lim_{n\to 0}\kappan(\alpha,\beta)=\kappa(\alpha,\beta)$.
\end{remark}

Let us denote by $\rho_n$ the probability measure defined with respect to $\xin$ in the same way as in \eqref{eq:rho_def}. We also set $\nu_n$ for the L\'evy measure of $\xin$ and $\hU_n$ for the renewal function of the ladder height process of $-\xin$.
The drift coefficient of $\Hn$ will be denoted by $a_n$.

\begin{lemma}\label{lem:rhon_weak_conv}
Under assumptions \ref{A:weak_conv} and \ref{A:EH_conv}, the sequence of probability measures $(\rho_n)$ converges weakly towards $\rho$ as $n \to \infty$.
\end{lemma}

\begin{proof}
Given the expression
\begin{equation*}
\rho_n(dx,dy)=\frac1{\EE \Hn_1}(\hU_n(x)\nu_n(x+dy)dx+a_n\delta_0(dx)\delta_0(dy))\,,
\end{equation*}
and assumption \ref{A:EH_conv}, it suffices to show that $\lim_{n\to\infty}a_n=a$ and that for all continuous functions $g\colon (0,\infty)^2\mapsto \R$ with compact support,
\begin{equation}\label{eq:31}
\lim_{n\to\infty}\int_{(0,\infty)^2}g(x,y)\hU_n(x)\nu_n(x+dy)dx=\int_{(0,\infty)^2}
g(x,y)\hU(x)\nu(x+dy)dx\,.
\end{equation}
Recall that from Lemma \ref{lem:WH_factors_conv}, $\lim_{n\to\infty}\kappan(0,\lambda)=\kappa(0,\lambda)$, for all $\lambda\geq 0$. Therefore, from \cite[Theorem 8.7]{Sato99}, $\lim_{n\to\infty}a_n=a$. Moreover, from \ref{A:weak_conv}
and the same theorem, we derive that for all bounded functions $h\colon (0,\infty)\mapsto\R$ with the support separated from the origin,
\begin{equation}\label{eq:nu_n_vague}
\lim_{n\to\infty}\int_{(0,\infty)}h(y)\,\nu_n(dy)=\int_{(0,\infty)}h(y)\,\nu(dy).
\end{equation}
Recall that for all $\lambda>0$,
\[
\int_{[0,\infty)}e^{-\lambda x}\,\hU_n(\od x)=\frac{1}{\hkappan(0,\lambda)}\quad
\text{and}\quad\int_{[0,\infty)}e^{-\lambda x}\,\hU(\od x)=\frac{1}{\hkappa(0,\lambda)},
\]
where $\hkappan$ and $\hkappa$ are the Laplace exponents of the upward ladder processes
of $-\xin$ and $-\xi$, respectively. In view of Lemma \ref{lem:WH_factors_conv} and the fact that $\widehat{H}$ is a subordinator that does not live on a lattice, $\hU(\od x)$ has no atom other than $0$ and it follows that
for all $x\geq0$,
\begin{equation}\label{eq:hU_n_conv}
\lim_{n\to\infty}\hU_n(x)=\hU(x).
\end{equation}
Let $g\colon(0,\infty)^2\mapsto\R$ be continuous with compact support included in
$[u,v]\times[c,d]$ with $0<u<v<\infty$ and $0<c<d<\infty$ and write
\[
\int_{(0,\infty)^2}g(x,y)\hU_n(x)\,\nu_n(x+\od y)\ud x=
\int_{(0,\infty)}\int_{(x,\infty)}g(x,z-x)\,\nu_n(\od z)\hU_n(x)\ud x.
\]
Since $z\mapsto g(x,z-x)$ is continuous with compact support in $(x,\infty)$, for all $x>0$ we have
\[
\lim_{n\to\infty}\int_{(x,\infty)}g(x,z-x)\,\nu_n(\od z)=
\int_{(x,\infty)}g(x,z-x)\,\nu(\od z),
\]
so that, for all $x>0$,
\[
\lim_{n\to\infty}\int_{(x,\infty)}g(x,z-x)\,\nu_n(\od z)\hU_n(x)=
\int_{(x,\infty)}g(x,z-x)\,\nu(\od z)\hU(x).
\]
Moreover,
\[
\int_{(x,\infty)}g(x,z-x)\,\nu_n(\od z)=0
\]
for $x\le u$ and $x\geq v$. Let $k$ be a bound for $g$.
From \eqref{eq:nu_n_vague} and \eqref{eq:hU_n_conv}, since $x\mapsto \nu[x,\infty)$ and $x\mapsto \hU(x)$ are
monotone, for all $x\in[u,v]$,
\begin{equation*}
    \int_{(x,\infty)}g(x,z-x)\,\nu_n(dz) \hU_n(x)\leq k\nu[u,\infty)\hU(v).
\end{equation*}
Then we derive \eqref{eq:31} from the dominated convergence.
\end{proof}


\begin{lemma}\label{2932}
Let $x>0$ and $(x_n)$ be any sequence such that $\lim_{n\rightarrow\infty}x_n=x$.
Assume that $\xin$ is a sequence of L\'evy processes that converges almost surely towards $\xi$.
Set $\tau^+_y=\inf\{t:\xi_t>y\}$ and $\tau^{+,n}_y=\inf\{t:\xin_t>y\}$, for any $y>0$. Then,
\[
\lim_{n\rightarrow\infty}\tau^{+,n}_{x_n}=\tau^+_x\,,\quad \text{a.s.}
\]
\end{lemma}
\begin{proof} Fix $x>0$ and let us first recall that almost surely, there is no $s<\tau^+_x$ such that
$\xi_s=x$ or $\xi_{s-}=x$, that is
\begin{equation}\label{3256}
\tau^+_x=\inf\{t:\xi_t\geq x\;\;\mbox{or}\;\;\xi_{t-}\geq x\},\quad \text{a.s.},
\end{equation}
see \cite[Lemma 3]{PecherskiiRogozin1969}. Now consider sequences $(\xin)$ and $(x_n)$ that satisfy the conditions of
the statement. Set 
$\oxi_t=\sup_{s\leq t}\xi_s$ and $\oxi_t^{(n)}=\sup_{s\leq t}\xi_s^{(n)}$, $t\geq0$. Recall that
for each $t\geq0$, $\xi_{t-}=\xi_t$, a.s.~and since $\oxi_t$ can only jump through a jump of $\xi$,
we derive that for all $t\geq0$, $\oxi_{t-}=\oxi_t$, a.s. Therefore, from  \cite[Theorem 6.2(i)]{1980_Whitt}, for all $t\geq0$,
\begin{equation}\label{3257}
\lim_{n\rightarrow\infty}\oxi_t^{(n)}=\oxi_t\,,\quad \text{a.s.}
\end{equation}
Let $t$ such that $t<\tau^+_x$. Then from \eqref{3256}, $\oxi_t<x$ so that from \eqref{3257},
for all $n$ large enough, $\oxi_t^{(n)}<x_n$. We derive that $t<\tau^{n,+}_{x_n}$, for all $n$ large
enough. Since $t$ is arbitrary, we obtain that $\liminf_{n\rightarrow\infty}\tau^{n,+}_{x_n}\geq\tau^+_x$,
a.s. Now let $t$ such that $t>\tau^+_x$. Then $\oxi_t>x$ and we derive from \eqref{3257} that
$\oxi_t^{(n)}>x_n$, for all $n$ large enough. This implies that $t\geq\tau^{+,n}_{x_n}$, for all $n$
large enough and since $t$ is arbitrary, we deduce that $\limsup_{n\rightarrow\infty}
\tau^{+,n}_{x_n}\leq \tau^+_x$, a.s.
\end{proof}

We denote by $\uPPn_x$, the law of
the process $\xin$ conditioned to stay positive and starting from $x\geq 0$. The following result is a
reinforcement of \cite[Theorem 2]{2005-LC-RD-ejp} and \cite[Theorem 3]{2018-Ivanovs-ann_appl_probab}.

\begin{lemma}\label{lem:cond_pos_uniform}
Let $x\geq 0$ and let any non negative sequence $(x_n)$ such that $\lim_{n\to\infty}x_n=x$.
When $x=0$ and $(x_n)$ is not identically zero from a certain rank onward, we assume that $0$ is regular
for $(0,\infty)$ for $\xi$ under $\PP$.

Then, under assumption \ref{A:weak_conv}, the family of probability measures $(\uPPn_{x_n})_{n\geq0}$
converges weakly toward $\uPP_x$ in the Skorokhod topology of $D$, as $n$ tends to $\infty$.
\end{lemma}

\begin{proof} Let us first consider the case $x>0$. Then we can assume that $x_n>0$ and $n\geq0$. We have to prove that for all $t>0$ and for all $\calF_t$-measurable, continuous, and bounded functional
$K$ in $D$,
\begin{equation}\label{eq:uEEn_conv}
\lim_{n\to\infty}\uEEn_{x_n}K(\xi)=\uEE_xK(\xi).
\end{equation}
Since $\hU_{n}$ and $\hU$ are increasing and continuous functions on $[0,\infty)$,
from Dini's lemma, the convergence \eqref{eq:hU_n_conv} is uniform on any compact set, which implies that
$\lim_{n\to\infty}\hU_{n}(x_n)=\hU(x)$. Let $c>0$. Appealing to Remark \ref{rem:representation_thm}, we can define, on some probability space, a family
of processes $\big(\xin\big)$ with distribution $\PPa{n}$ that converges almost surely towards $\xi$ with distribution $\PP$ on $D$. Then it follows from Lemma \ref{2932} applied to the processes $(-\xin)$ and $-\xi$ that  $\lim_{n\rightarrow\infty}\tau_{-x_n}^{n,-}=\tau_{-x}^-$, a.s., where $\tau_{-x_n}^{n,-}=\inf\{t:\xi_t^{(n)}<-x_n\}$. We derive from this almost sure convergence, the uniform convergence of $(\hU_{n})$ towards $\hU$ on compact sets and the continuity of $K$ that almost surely,
\[\lim_{n\to\infty}\hU_{n}(\xi_t^{(n)}+x_n)K(\xin+x_n)
\ind_{\{\xi_t^{(n)}+x_n\le c\}}\ind_{\{t<\tau_{-x_n}^{n,-}\}}=\hU(\xi_t+x)K(\xi+x)
\ind_{\{\xi_t+x\le c\}}\ind_{\{t<\tau_{-x}^-\}}\,,\]
so that from dominated convergence,
\[
\lim_{n\to\infty}\EEa{n}_{x_n}\Big[\hU_{n}(\xi_t)
K(\xi),\xi_t\leq c,t<\tau_{0}^-\Big]=\EE_x\Big[\hU(\xi_t)K(\xi),\xi_t\le c,t<\tau_{0}^-\Big],
\]
where $\PP^{(n)}_x$ (with expectation $\EE^{(n)}_x$) denotes the probability under which $\xi^{(n)}_0=x$, almost surely. Then we have proved that,
\[
\lim_{n\to\infty}\uEEn_{x_n}\big[K(\xi),\xi_t\leq c\big]=
\uEE_x\big[K(\xi),\xi_t\le c\big].
\]
Now recall, c.f. \cite[Lemma 1]{2005-LC-RD-ejp}, that
\[
\EEa{n}_{x_n}\Big[\hU_{n}(\xi_t),\xi_t\leq c,t<\tau_{0}^-\Big]+
\EEa{n}_{x_n}\Big[\hU_{n}(\xi_t),\xi_t>c,t<\tau_{0}^-\Big]=\hU_{n}(x_n),
\]
and that
\[
\lim_{n\to\infty}\EEa{n}_{x_n}\Big[\hU_{n}(\xi_t),\xi_t\leq c,t<\tau_{0}^-\Big]=
\EE_x\Big[\hU(\xi_t),\xi_t\leq c,t<\tau_{0}^-\Big],
\]
so that
\[
\lim_{n\to\infty}\EEa{n}_{x_n}\Big[\hU_{n}(\xi_t),\xi_t>c,t<\tau_{0}^-\Big]=
\hU(x)-\EE_x\Big[\hU(\xi_t),\xi_t\le c,t<\tau_{0}^-\Big]=\EE_x\Big[\hU(\xi_t),\xi_t>c,t<\tau_{0}^-\Big].
\]
Then we have proved that
\[
\lim_{n\to\infty}\frac1{\hU_{n}(x_n)}
\EEa{n}_{x_n}\Big[\hU_{n}(\xi_t),\xi_t>c,t<\tau_{0}^-\Big]=
\lim_{n\to\infty}\uPPn_x(\xi_t>c)=\uPP_x(\xi_t>c).
\]
Let $\eta>0$ and $N$ be a bound for $K$. Then we derive from the above limit that for $c$ sufficiently large,
\[
\limsup_{n\to\infty}\uEEn_{x_n}[K,\xi_t>c]\leq N\uPP_x(\xi_t>c)<\eta.
\]
Therefore, we can conclude that \eqref{eq:uEEn_conv} holds since we have already proved that
$\lim_{n\to\infty}\uEEn_{x_n}[K,\xi_t\leq c]=\uEE_x[H,\xi_t\leq c]$.

Let us now consider the case $x=0$. Let $y>0$ and recall that for all $n\geq0$, $\uPPn_0$ is
the law of the post-infimum process under $\uPPn_y$. Since the time
$m=\inf\{t:\xi_t=\inf_{s\geq0}\xi_s\}$, where the infimum is attained, is a.s.~unique under $\uPP_y$,
we derive from the weak convergence of $(\uPPn_{y})_{n \geq 0}$ towards $\uPP_{y}$ proved
above and the continuous mapping theorem that the post-infimum process under $\uPPn_y$ converges
weakly towards the post-infimum process under $\uPP_y$ (it follows from the same arguments as in the last paragraph  of the proof of Lemma \ref{lem:WH_factors_conv}). As a consequence, the sequence of probability
measures $(\uPPn_0)$ converges weakly toward $\uPP_0$, as $n$ tends to $\infty$ and
this proves our result in the case where $x=0$ and $x_n=0$ from a certain rank.\\

Now assume that $(x_n)$ is not identically zero from a certain rank onwards. Then from what has just been proved,
without loss of generality we can assume that $x_n>0$ for all $n$. Then assume also that $0$ is regular for
$(0,\infty)$ for $\xi$ under $\PP$, so that $\uPP_{0}(\xi_0=0)=1$. Let us consider some
probability space $(\Omega,{\calF},P)$, on which we can define a sequence $(\xi^{x_n,n,\uparrow})$ of processes
such that for each $n$, $\xi^{x_n,n, \uparrow}$ has law $\uPPn_{x_n}$. Let also
$(\xi^{0,n,\uparrow})$ be a sequence of processes that is independent of the sequence $(\xi^{x_n,n,\uparrow})$ and
such that $\xi^{0,n,\uparrow}$ has law $\uPPn_{0}$ for each $n$. From what has been proved above and
Skorokhod's representation theorem, we can assume that $(\xi^{0,n,\uparrow})$ converges almost surely on $D$ toward
a process with law $\uPP_{0}$. It follows from \cite[Theorem 1]{2005-LC-RD-ejp} that the process $\xi^{(n),\uparrow}$
defined for each $t\geq 0$, by
\[
\xi^{(n),\uparrow}_t=
\begin{cases}
\xi^{x_n,n,\uparrow}_t,&t<m_{x_n}^{(n)},\\
\xi^{0,n,\uparrow}_{t-m_{x_n}^{(n)}}+\xi^{x_n,n,\uparrow}_{m_{x_n}^{(n)}}\wedge\xi^{x_n,n,\uparrow}_{m_{x_n}^{(n)}-},
&t\geq m_{x_n}^{(n)},
\end{cases}
\]
where $m_{x_n}^{(n)}$ is the time at which $\xi^{x_n,n,\uparrow}$ attains its overall infimum, has law
$\uPPn_{x_n}$, see also \cite[proof of Theorem 2]{2005-LC-RD-ejp}. From this construction and
since $(\xi^{0,n,\uparrow})$ converges a.s.~towards a process with law $\uPP_{0}$, the sequence of processes
$(\xi^{(n),\uparrow})$ will also converge weakly to $\uPP_{0}$ if both sequences $(m_{x_n}^{(n)})$ and
$(M_n)$, where $M_n:=\sup_{t\leq m_{x_n}^{(n)}}\xi^{x_n,n,\uparrow}_t$ converge in probability to $0$ as $n$ tends to $\infty$. But we can prove in the same way as in \cite[proof of Theorem 2]{2005-LC-RD-ejp} that, for all $\eta>0$,
\[
P(m_{x_n}^{(n)}>\eta)\le \frac{\hU_n(x_n)-\hU_n(0)}{\hU_n(x_n)}P\left(\inf_{t\leq \eta}\xin_t\geq 0\right)
\quad\text{and}\quad P(M_n>\eta)\leq \frac{\hU_n(x_n)-\hU_n(0)}{\hU_n(\eta)}\,,
\]
where $\xin$ has law $\PPa{n}$ under $P$. From uniform convergence, $(\hU_n(x_n)-\hU_n(0))$ converges
to 0, so that $(M_n)$ tends to $0$ in probability. If $0$ is not regular for $(-\infty,0)$ under $\PP$, then
$\hU_n(x_n)$ tends to $\hU(0)>0$, so that $\left(\frac{\hU_n(x_n)-\hU_n(0)}{\hU_n(x_n)}\right)$ tends to 0. If $0$ is regular for $(-\infty,0)$ under $\PP$, then $\left(P\left(\inf_{t\leq \eta}\xin_t\geq 0\right)\right)$
tends to $0$ and  $\frac{\hU_n(x_n)-\hU_n(0)}{\hU_n(x_n)}\leq 1$ for all $n$. Hence, in both cases,
$\big(m_{x_n}^{(n)}\big)$ tends to $0$ in probability, as required.
\end{proof}

Let $\txin$ be the process that is constructed from $\xin$ in the same way as
$\txi$ is constructed from $\xi$, see \eqref{1258}.

\begin{lemma}\label{lem:oxi_conv_in_oD}
Under the assumptions \ref{A:weak_conv} and \ref{A:EH_conv}, the sequence of processes $(\txin)_{n\geq0}$ converges weakly in
$\tD$ towards the process $\txi$, as $n$ tends to $\infty$.
\end{lemma}

\begin{proof}
Let us recall the bijection $\varphi:\tomega\mapsto (\omega,\omega')$ from $\tD$ to
$D\times D$ given in \eqref{bijection} by
\[\tomega(t)=\left\{\begin{array}{ll}
\omega'(t)&\mbox{if $t\geq 0$,}\\
-\omega(-t-)&\mbox{if $t<0$.}
\end{array}\right.\]
Let us also recall from \eqref{1258} the definition of the law $\calP$ of the process $\txi$ and let $\calP_n$ denote the law of $\txin$. Then from the construction described above, these distributions are respectively defined by
\begin{align*}
\varphi\circ\calP_n(\od\omega,\od\omega')&=
\int_{[0,\infty)^2}\rho_n(\od x,\od y)\,\uPPn_x(\od\omega)\,\PP_y^{(n)}(\od\omega'),\intertext{and}
\varphi\circ\calP(\od\omega,\od\omega')&=
\int_{[0,\infty)^2}\rho(\od x,\od y)\,\uPP_x(\od\omega)\,\PP_y(\od\omega'),
\end{align*}
see Subsection 2.4 in \cite{2011-JB-MS-bulllms}. Note that $\varphi$ is a continuous bijection. Then 
the weak convergence of $\calP_n$ towards $\calP$ is equivalent to 
the weak convergence of $\varphi\circ\calP_n$ towards $\varphi\circ\calP$. In order to prove the latter convergence, it is enough to prove that for any bounded and continuous functional $H\colon D\times D\rightarrow\mathbb{R}$ of the form $H(\omega,\omega')=H_1(\omega)H_2(\omega')$, where 
$H_i\colon D\rightarrow\mathbb{R}$, $i=1,2$ are bounded and continuous functionals, 
\begin{equation}\label{6320}
\lim_{n\to\infty}\varphi\circ\calE_n(H)=
\lim_{n\to\infty}\calE_n(H\circ \varphi)=\calE(H\circ \varphi),
\end{equation}
where $\calE_n$ (resp. $\calE$) denotes the expectation associated to $\calP_n$ (resp. $\calP$). By definition,
\[\varphi\circ\calE_n(H)=\varphi\circ\calE_n(H_1H_2)=
\int_{[0,\infty)^2}\rho_n(\od x,\od y)\uEEn_x [H_1]\EE_y^{(n)}[H_2].\]
Then fix $H_1$ and $H_2$ as above and define the sequence of functions
$f_n\colon [0,\infty)^2\rightarrow\mathbb{R}$ by 
$f_n(x,y)=\uEEn_x [H_1]\EE_y^{(n)}[H_2]$.
It follows from the weak convergence of $\xin$ towards $\xi$, as $n$ tends to $\infty$ that for
all $y\geq0$ and all sequences $(y_n)$ that converge towards $y$ that
$\lim_{n\to\infty}\EE_{y_n}^{(n)}[H_2]=\EE_y[H_2]$. Then we derive from Lemma \ref{lem:cond_pos_uniform} that for all $(x,y)\in(0,\infty)\times[0,\infty)$ and all sequences $(x_n,y_n)\in[0,\infty)^2$, $n\geq 0$ that converge toward $(x,y)$, 
\begin{equation}\label{2994}
\lim_{n\rightarrow\infty}f_n(x_n,y_n)=f(x,y)\coloneqq \uEE_x [H_1]\EE_y[H_2],
\end{equation}
and the result is still true for all 
$(x,y)\in\{0\}\times[0,\infty)$, 
provided that when $x_n>0$ for infinitely many indices 0 is regular for $(0,\infty)$, for $\xi$.
From Lemma \ref{lem:rhon_weak_conv} and Skorokhod's representation theorem, there is a probability space $(\Omega',\mathcal{F}',P)$ on which we can define random variables $\Lambda_n$, $n\geq 0$ and $\Lambda$
with values in $[0,\infty)^2$ and respective distributions 
$\rho_n$, $n\geq 0$ and $\rho$ such that the sequence $(\Lambda_n)$
converges almost surely towards $\Lambda$. Moreover, note that from 
\eqref{eq:rho_def}, $\rho(\{0\}\times[0,\infty))=a/{\bf E}(H_1)$, where $a$ is the drift of $H$. If 0 is not regular for $(0,\infty)$, for $\xi$,  then $H$ is a compound Poisson process and hence $a=0$. So we can claim from \eqref{2994} that in any case, 
$\lim_{n\rightarrow\infty}f_n(\Lambda_n)=f(\Lambda)$, almost surely. 
Recall that the functions $f_n$, $n\geq 0$ and $f$ are continuous and bounded (uniformly in $n$). Then it follows from dominated convergence that $\lim_{n\rightarrow\infty}E(f(\Lambda_n))=E(f(\Lambda))$, where $E$ denotes the expectation with respect to $P$. This provides \eqref{6320} since $E(f_n(\Lambda_n))=\varphi\circ\calE_n(H)$ and $E(f(\Lambda))=\varphi\circ\calE(H)$. 
\end{proof}

\begin{lemma}\label{lem:oAn_conv}
Let us define $\tA^{(n)}_t=\int_{-\infty}^te^{\txin_s}\,ds$, $n\geq0$ and
$\tA_t=\int_{-\infty}^te^{\txi_s}\,ds$, $t\in\R$.
Then under assumptions \ref{A:weak_conv} and \ref{A:EH_conv}, the sequence of processes $(\tA^{(n)}_t)_{t\in\R}$,
$n\geq0$ converges weakly in $\tD$ towards the process $(\tA_t)_{t\in\R}$, as $n$ tends
to $\infty$.
\end{lemma}

\begin{proof}
For any $t,t_0\in\mathbb{R}$ such that $t\geq t_0$, we set $\tA_{t_0,t}=\int_{t_0}^te^{\txi_s}\ud s$ and
$\tA_{t_0,t}^{(n)}=\int_{t_0}^te^{\txin_s}\ud s$.
From Lemma \ref{lem:oxi_conv_in_oD} and the Skorokhod representation theorem applied to $\tD$ (see Remark \ref{rem:representation_thm}), one can construct a family of processes
$(\txin)_{n\geq0}$ (on some probability space) such that $e^{\txin}$ converges a.s.
towards $e^{\txi}$, as $n$ tends to $\infty$. Note that the function
$\omega\mapsto\int_{t_0}^t\omega_u\,du$ is continuous on $D([t_0,t])$.
Then  $\tA_{t_0,t}^{(n)}$ converges a.s. towards $\tA_{t_0,t}$ as $n$ tends to $\infty$.
On the other hand, from \cite[proof of part $(i)$ of Corollary 4]{2011-JB-MS-bulllms},
\begin{equation*}
\E\left[\tA^{(n)}_0\right]=
\frac{1}{\EE H_1^{(n)}}\int_0^{\infty}e^{-y}\hU_n(y)\ud y=
\frac{1}{\EE H_1^{(n)}\hkappan(0,1)}.
\end{equation*}
Therefore, from assumptions \ref{A:weak_conv} and \ref{A:EH_conv} and Lemma \ref{lem:WH_factors_conv}, $\lim_{n\to\infty}\E \Big[\tA^{(n)}_0\Big]=
\E \Big[\tA_0\Big]$. Moreover, from the almost sure convergence assumed above
and dominated convergence, for all $t_0\leq 0$,
$\lim_{n\to\infty}\E\Big[\tA_{t_0,0}^{(n)}\Big]=
\E\Big[\tA_{t_0,0}\Big]$. It follows that for all $t_0\leq 0$,
$\lim_{n\to\infty}\E\Big[\tA_{t_0}^{(n)}\Big]=
\E\Big[\tA_{t_0}\Big]$.

Let $k\geq 1$, $t_i\in\mathbb{R}$, $i=0,\ldots,k$ such that $-\infty<t_1\leq \ldots\leq t_k<\infty$ and $t_0\leq\min(0,t_1)$. Let also $\delta\in(0,1)$. Then from the inequality
$1-e^{-x}\leq \big(1-e^{-(1-\alpha)x/(1-\delta)}\big)+\big(1-e^{-\alpha x/\delta}\big)$, $x\geq 0$, $\alpha\in(0,1)$ (obtained by integrating $e^{-y}\ind_{\{y<x\}}\leq e^{-y}\ind_{\{(1-\delta)y<x(1-\alpha)\}}+e^{-y}\ind_{\{\delta y<\alpha x\}}$) we derive that
for all $\lambda_1,\dots,\lambda_k\geq 0$,
\begin{align*}
\EE \left[(1-e^{-\lambda_1 \tA_{t_1}^{(n)}})
e^{-\lambda_2\tA_{t_1,t_2}^{(n)}-\ldots-\lambda_k\tA_{t_{k-1},t_k}^{(n)}}\right]&\leq\EE\left[(1-e^{-\lambda_1\tA_{t_0,t_1}^{(n)}/(1-\delta)})e^{-\lambda_2\tA_{t_1,t_2}^{(n)}-\ldots-\lambda_k\tA_{t_{k-1},t_k}^{(n)}}\right]\\
&+\EE\left[1-e^{-\lambda_1\tA_{t_0}^{(n)}/\delta}\right].
\end{align*}
Moreover, from what has been proved above,
\[
\EE\left[1-e^{-\lambda_1\tA_{t_0}^{(n)}/\delta}\right]\leq \frac{\lambda_1}{\delta}\EE\Big[\tA_{t_0}^{(n)}\Big]\longrightarrow \frac{\lambda_1}{\delta}\EE\Big[\tA_{t_0}\Big] \quad \text{as} \quad n \to \infty.\]
Note also that from the dominated convergence, $\lim_{t_0\to-\infty}\EE\big[\tA_{t_0}\big]=0$.
Therefore, for all $\varepsilon>0$ and $\delta\in(0,1)$,
there is $t_0\leq0$ and an integer $n_0\geq0$ such that for all $t\leq t_0$ and $n\geq n_0$,
\begin{align*}
&\EE\left[\left(1-e^{-\lambda_1 \tA_{t_1}^{(n)}}\right)
e^{-\lambda_2\tA_{t_1,t_2}^{(n)}-\ldots-\lambda_k\tA_{t_{k-1},t_k}^{(n)}}\right]
\\ & \qquad \leq  \EE\left[\left(1-e^{-\lambda_1\tA_{t,t_1}^{(n)}/(1-\delta)}\right)e^{-\lambda_2\tA_{t_1,t_2}^{(n)}-\ldots-\lambda_k\tA_{t_{k-1},t_k}^{(n)}}\right]+\epsilon.\end{align*}
Now recall that from our assumption, $\Big(\tA_{t,t_1}^{(n)},\tA_{t_1,t_2}^{(n)},\dots,\tA_{t_{k-1},t_k}^{(n)}\Big)$ converges almost surely towards $\big(\tA_{t,t_1},\tA_{t_1,t_2},\ldots,\tA_{t_{k-1},t_k}\big)$. Then we derive from the previous inequality that for all $t\leq t_0$,
\begin{align*}
\limsup_{n\to\infty}&\EE\left[\left(1-e^{-\lambda_1 \tA_{t_1}^{(n)}}\right)
e^{-\lambda_2\tA_{t_1,t_2}^{(n)}-\ldots-\lambda_k\tA_{t_{k-1},t_k}^{(n)}}\right]\\ \leq
&\EE\left[\left(1-e^{-\lambda_1\tA_{t,t_1}/(1-\delta)}\right)e^{-\lambda_2\tA_{t_1,t_2}-\ldots-\lambda_k\tA_{t_{k-1},t_k}}\right]+\varepsilon.
\end{align*}
Note that for all $t\leq t_1$, $\tA_{t,t_1}\leq \tA_{t_1}$. Then by taking $\varepsilon$ and $\delta$ to 0, it follows that,
\begin{equation}\label{7399}
\limsup_{n\to\infty}\EE\left[\left(1-e^{-\lambda_1 \tA_{t_1}^{(n)}}\right)e^{-\lambda_2\tA_{t_1,t_2}^{(n)}-\ldots-\lambda_k\tA_{t_{k-1},t_k}^{(n)}}\right]\leq
\EE\left[\left(1-e^{-\lambda_1\tA_{t_1}}\right)e^{-\lambda_2\tA_{t_1,t_2}-\ldots-\lambda_k\tA_{t_{k-1},t_k}}\right].
\end{equation}
On the other hand, since $\tA_{t_1}^{(n)}\geq\tA_{t,t_1}^{(n)}$, we have for all $n$ and $t\leq t_1$,
\[
\EE\left[\left(1-e^{-\lambda_1\tA_{t,t_1}^{(n)}}\right)e^{-\lambda_2\tA_{t_1,t_2}^{(n)}-\ldots-\lambda_k\tA_{t_{k-1},t_k}^{(n)}}\right]\leq \EE\left[\left(1-e^{-\lambda_1\tA_{t_1}^{(n)}}\right)e^{-\lambda_2\tA_{t_1,t_2}^{(n)}-\ldots-\lambda_k\tA_{t_{k-1},t_k}^{(n)}}\right],
\]
which, thanks to the almost sure convergence of
$\Big(\tA_{t,t_1}^{(n)},\tA_{t_1,t_2}^{(n)},\ldots,\tA_{t_{k-1},t_k}^{(n)}\Big)$ towards \linebreak $\big(\tA_{t,t_1},\tA_{t_1,t_2},\ldots,\tA_{t_{k-1},t_k}\big)$ leads to
\[
\EE\left[\left(1-e^{-\lambda_1\tA_{t,t_1}}\right)e^{-\lambda_2\tA_{t_1,t_2}-\ldots-\lambda_k\tA_{t_{k-1},t_k}}\right]\leq \liminf_{n\to\infty}\EE\left[\left(1-e^{-\lambda_1\tA_{t_1}^{(n)}}\right)e^{-\lambda_2\tA_{t_1,t_2}^{(n)}-\ldots-\lambda_k\tA_{t_{k-1},t_k}^{(n)}}\right].
\]
Taking $t$ to $-\infty$ we obtain,
\begin{equation}\label{7339}
\EE\left[\left(1-e^{-\lambda_1\tA_{t_1}}\right)e^{-\lambda_2\tA_{t_1,t_2}-\ldots-\lambda_k\tA_{t_{k-1},t_k}}\right]\leq \liminf_{n\rightarrow\infty}\EE\left[\left(1-e^{-\lambda_1\tA_{t_1}^{(n)}}\right)e^{-\lambda_2\tA_{t_1,t_2}^{(n)}-\ldots-\lambda_k\tA_{t_{k-1},t_k}^{(n)}}\right].
\end{equation}
From \eqref{7399} and \eqref{7339} we derive that for all $k\geq 1$ and $\lambda_1,\dots,\lambda_k\geq 0$,
\[
\lim_{n\to\infty}\EE\left[e^{-\lambda_1\tA_{t_1}^{(n)}}e^{-\lambda_2\tA_{t_1,t_2}^{(n)}-\ldots-\lambda_k\tA_{t_{k-1},t_k}^{(n)}}\right]=\EE\left[e^{-\lambda_1\tA_{t_1}}e^{-\lambda_2\tA_{t_1,t_2}-\ldots-\lambda_k\tA_{t_{k-1},t_k}}\right].
\]
This proves that the process $(\tA_{t}^{(n)})_{t\in\R}$ converges toward the process $(\tA_{t})_{t\in\R}$,
in the sense of finite dimensional distributions, as $n$ tends to $\infty$. Since the processes
$(\tA_{t})_{t\in\R}$ and $(\tA_{t}^{(n)})_{t\in\R}$, $n\geq0$
are increasing and continuous, the convergence holds weakly in the topology of $\tD$, see for instance \cite[Corollary 1.3]{89-Aldous-ann_probab} applied to $B^{(n)}_t:=\tA_{\log t}^{(n)}$.
\end{proof}
We are now ready to prove Theorem \ref{thm:ssMps_0_conv}.
\begin{proof}[Proof of Theorem \ref{thm:ssMps_0_conv}]
Assume first that \ref{A:weak_conv} and \ref{A:EH_conv} hold. From \cite[Theorem 7.2]{1980_Whitt} (set $B^{(n)}_t\coloneq \tA_{\log t}^{(n)}$)
and an application of Lemma \ref{lem:oAn_conv}, we obtain that the processes
$\sigma_t^{(n)}=\inf\{s\in\R:\tA_s^{(n)}>t\}$, $t\geq0$ converge weakly in $D([0,\infty),\R)$
toward the process $\sigma_t=\inf\{s \in \R:\tA_s>t\}$, $t\geq0$, as $n$ tends to $\infty$.
Then we derive that $(X^{(n)})$ converges weakly toward the process $X$ from the transformation \eqref{eq:pssMp_ex} applied to $\txin$ and $\txi$, Lemma \ref{lem:oxi_conv_in_oD} and \cite[Theorem 3.1]{1980_Whitt}.

Assume now that the sequence of processes $(X^{(n)})$ converges weakly in the Skorokhod topology of $D$ towards the process $X$ as $n$ tends to $\infty$. Set $\gamma=\inf\{t:X_t\geq1\}$ and
$\gamma_n=\inf\{t: X^{(n)}_t\geq1\}$. Then note that from \eqref{eq:pssMp_ex},
\[
\sigma_\gamma=\inf\{t:\txi_t\geq0\}=0 \quad\text{a.s.}\quad \text{and}\quad\sigma^{(n)}_{\gamma_n}=
\inf\{t:\txi_t^{(n)}\geq0\}=0\quad \text{a.s.}
\]
Moreover, for all $s,t\in\R_+$,
\[
\sigma_t-\sigma_s=\int_s^t\big(X_u\big)^{-1}\ud u\quad \text{and}\quad
\sigma_t^{(n)}-\sigma_s^{(n)}=\int_s^t\Big(X_u^{(n)}\Big)^{-1}\ud u,
\]
so that for all $t\geq0$,
\begin{equation}\label{eq:sigma_gamma}
\sigma_t=\int_\gamma^t\big(X_u\big)^{-1}\ud u\quad \text{and}\quad \sigma_t^{(n)}=\int_{\gamma_n}^t\Big(X_u^{(n)}\Big)^{-1}\ud u.
\end{equation}
We derive from our assumption of weak convergence and \cite[Proposition VI.2.11, p. 341]{JS87} that $({\gamma_n},X^{(n)})$ converges weakly toward $(\gamma, X)$, as $n$ tends to $\infty$. Since the function $\omega\mapsto\int_{s}^t\omega_u\ud u$ is continuous on $D([s,t])$, for all $0\leq s\leq t$, we derive from \eqref{eq:sigma_gamma} and weak convergence of $\big({\gamma_n},X^{(n)}\big)$ that the increasing processes $\sigma^{(n)}$, $n\geq0$ converge weakly towards $\sigma$, as $n$ tends to $\infty$ (here again, we can assume almost sure convergence via Skorokhod's representation theorem in order to justify this convergence). Then we derive from \eqref{eq:pssMp_ex} and \cite[Theorems 7.2 and 5.1(i)]{1980_Whitt} that the processes $\txin$, $n\geq0$, converge weakly towards $\txi$ as $n \to \infty$, so that \ref{A:weak_conv} is satisfied.

Next, note that $X_\gamma=e^{\txi_0}$ and $X^{(n)}_\gamma=e^{\txin_0}$.
We can check from \cite{2011-JB-MS-bulllms} (see the definition of $\rho_1$ on p.~101) that
\begin{align*}
\PP\Big(e^{\txi_0}\in[a,b]\Big)&=
\frac{1}{\EE H_1}\int_a^b\hU(x)\nu[x,\infty)\ud x\intertext{and}
\PP\Big(e^{\txin_0}\in[a,b]\Big)&=\frac{1}{\EE H_1^{(n)}}\int_a^b\hU_n(x)\nu_n[x,\infty)\ud x,
\end{align*}
for all $0<a<b<\infty$, and therefore,
\[
\lim_{n\to\infty}\frac{1}{\EE H_1^{(n)}}\int_a^b\hU_n(x)\nu_n[x,\infty)\ud x=
\frac{1}{\EE H_1}\int_a^b\hU(x)\nu[x,\infty)\ud x.
\]
On the other hand, since \ref{A:weak_conv} is satisfied, we know from \eqref{eq:nu_n_vague}, \eqref{eq:hU_n_conv} and monotonicity of $x \mapsto \hU(x)$ and $x \mapsto \nu([x,\infty))$ that
\[
\lim_{n\to\infty}\int_a^b\hU_n(x)\nu_n([x,\infty))\ud x=\int_a^b\hU(x)\nu([x,\infty))\ud x.
\]
This proves \ref{A:EH_conv}.
\end{proof}
\subsection{Discussion of assumptions}
Now let us discuss in details assumptions of Theorem \ref{thm:ssMps_0_conv}. While \ref{A:weak_conv} is rather obvious and expected for Theorem \ref{thm:ssMps_0_conv}, the convergence of means of the corresponding ladder height processes depicted in \ref{A:EH_conv} is much more elusive. For this reason, let us first construct a counterexample to illustrate its necessity:
\begin{example}\label{ex:counterex_ssMps_conv_0}
    Consider a sequence of driftless subordinators $\big(S^{(n)}\big)_{n \in \N}$ defined by their L\'{e}vy measures $\nu_n(\od x) = \nu_n(x)\ud x$ as follows:
    \[
    \nu_n(x)=\begin{cases}
        \frac1x, & x \in (0,1),\\
        n^{-3/2}, & x \in [1,n),\\
        0, & x \geq n,
    \end{cases}
    \]
    where we accept a slight ambiguity by using the same symbol $\nu_n$ to denote L\'{e}vy measure and its density. By \cite[Theorem 8.7]{Sato99} and \cite[Corollary VII.3.6]{JS87} we see that $S^{(n)} \to S$ in the Skorokhod space as $n \to \infty$, where $S$ is given by its L\'{e}vy density $\nu(\od x)=\nu(x)\ud x$ and
    \[
    \nu(x) = \frac1x \ind_{(0,1)}(x).
    \]
    Then one may verify that $\EE S_1=1$ and
    \[
    \EE S_1^{(n)} = 1+\frac{n^2-1}{2n^{3/2}} \to \infty \quad \text{as} \quad n \to \infty,
    \]
    so that $\EE S_1^{(n)}$ does not converge to $\EE S_1$ as $n \to\infty$.

    Next, since the densities of $S^{(n)}$ are non-increasing, Vigon's theory of friendship asserts that these subordinators are upward ladder height processes of spectrally positive L\'{e}vy processes $Z^{(n)}$ and $Z$ such that $Z^{(n)}$ converges in the Skorokhod space to $Z$ as $n \to \infty$, see the first proposition of \cite[Subsection 7.3]{Vigon02_PhD}. Thus, \ref{A:EH_conv} does not hold for $\big(Z^{(n)}\big)_{n \in \N}$ and $Z$. This example shows that Theorem \ref{thm:ssMps_0_conv} does not hold without \ref{A:EH_conv}.
\end{example}
Let us further explore Example \ref{ex:counterex_ssMps_conv_0}. First note that all $\nu_n$ and $\nu$ have bounded supports, so that all exponential moments of $S^{(n)}$ and $S$ are well defined, and their Laplace exponents are analytic on $\C$ (see \cite[Lemmas 25.6 and 25.7, Theorem 25.17]{Sato99}). Since these are Laplace exponents of upward ladder height processes of $Z^{(n)}$ and $Z$, let us denote them by $\phin$ and $\phi$, respectively. Then one may derive that
    \[
    \phi(\lambda) = c_\lambda \quad \text{and} \quad \phin(\lambda) = \begin{cases} c_\lambda - \frac1\lambda\cdot\frac{(e^{-\lambda}-e^{-\lambda n})}{n^{3/2}} -\frac{n-1}{n^{3/2}}, & \lambda \neq 0, \\ 0,& \lambda=0,\end{cases}
    \]
    where $c_\lambda = \int_0^1 (1-e^{-\lambda x})\,\frac{\od x}{x}$. Thus, we see that $\phin(\lambda) \to \phi(\lambda)$ if and only if $\lambda \geq 0$, and otherwise $\phin(\lambda) \to -\infty$ as $n \to \infty$. In terms of exponential moments, we get
    \[
    \lim_{n \to \infty} \EE e^{\beta S^{(n)}_1} = \EE e^{\beta S_1} =e^{c_\beta}\quad \text{if and only if} \quad \beta<0,
    \]
    and otherwise, for all $\beta>0$,
    \[
    \lim_{n \to \infty} \EE e^{\beta S^{(n)}_1} = \infty.
    \]
    Note that the descending ladder height subordinators are all deterministic drifts, and so $\hphin(z)=\hphi(z)=z$ for all $z \in \C$. It follows that $Z^{(n)}$ and $Z$ also have all exponential moments, and denoting by $\psin$ and $\psi$ their Laplace exponents, i.e.
\[
\EE e^{z \xin_t} = e^{t\psin(z)} \quad \text{and} \quad \EE e^{z \xi_t} = e^{t\psi(z)}, \quad z \in \C,
\]
the Wiener-Hopf factorisation reads as follows 
\[
\psin(z) = -z\phin(-z) \quad \text{and} \quad \psi(z)=-z\phi(-z), \quad z \in \C,
\]
see e.g. \cite[Th\'{e}or\`{e}me 4.2.2]{Vigon02_PhD}.
Hence we conclude that $\psin(\beta) \to \psi(\beta)$ if and only if $\beta \leq 0$ and consequently,
\[
\lim_{n \to \infty}\EE e^{\beta \xin_1} = \EE e^{\beta \xi_1} = e^{\beta c_\beta} \quad \text{if and only if} \quad \beta \leq 0,
\]
and otherwise, for all $\beta>0$,
\[
    \lim_{n \to \infty} \EE e^{\beta \xin_1} = \infty.
\]
We also see that $\EE \xin_1=\EE \xi_1=0$.

As we have seen in the above example, the finiteness of $\EE \Hn_1$, and in fact, even the finiteness of $\EE e^{\beta \Hn_1}$ for all $\beta \in \R$, is not enough for \ref{A:EH_conv} to hold. In view of \cite[Theorem 3.6]{Billingsley99}, the sufficient condition is the uniform integrability of $\Hn_1$, but it would be convenient to have some elegant condition in terms of $\xin$ only, since the ladder height processes are usually much less tractable and explicit formulas are rarely available. In Proposition \ref{prop:pointwise} we provide one such condition which will be useful later on. Recall that
\[
\C_a = \{z \in \C \colon \Re z=a\}, \quad \C_{(a,b)} = \{z \in \C \colon a < \Re z < b\}, \quad \C_{[a,b]}=\rA_{(a,b)} \cup \C_a \cup \C_b,
\]
and
\[
\rAC = \{ z \in \C : \EE |e^{z\xi_1}|<\infty\}.
\]
It is well known that $\rAC$ is a complex vertical strip containing $i\R$, c.f. e.g. \cite[Theorem 25.17]{Sato99}.
\begin{proposition}\label{prop:pointwise}
    Assume a sequence of L\'{e}vy processes $(\xin)_{n \in \N}$ converging in the Skorokhod space $n \to \infty$ to a L\'{e}vy process $\xi$ which is not a compound Poisson process. Suppose that $\limsup_{t \to \infty} \xi_t=+\infty$ a.s. or $\lim_{t \to \infty} \xin_t=-\infty$ a.s. from a certain rank. Furthermore, assume that there is $\delta>0$ such that $\EE e^{\delta \xi_1}<\infty$, $\EE e^{\delta \xin_1}<\infty$ for all $n \in \N$ and $\psin(\lambda) \to \psi(\lambda)$ as $n \to \infty$ for all $\lambda \in [0,\delta]$. Then the family $\Big(\Hn_1\Big)_{n \in \N}$ is uniformly integrable and
    \[
    \lim_{n \to \infty} \EE \Hn_1=\EE H_1 \in (0,\infty).
    \]
\end{proposition}
\begin{proof}
    Our assumptions imply in particular that $\C_{[0,\delta]} \subseteq \rAC^{(n)}$ and $\C_{[0,\delta]} \subseteq \rAC$. The analytic extension of the Wiener-Hopf factorisation (see e.g. \cite[Th\'{e}or\`{e}me 4.2.2]{Vigon02_PhD} gives
    \[
    \psin(z) = -\phin(-z)\hphin(z), \quad z \in \C_{[0,\delta]}
    \]
    and according to \cite[Lemma A.1]{PS18-ejp} one may take $\phin(z)=\kappan(0,z)$ and $\hphin(z)=\hkappan(0,z)$ (note that $\xi$ is not a compound Poisson process, so $h(q)=1$ for all $q \geq 0$).
    Now, Lemma \ref{lem:WH_factors_conv} directly implies $\lim_{n \to \infty} \phin(\lambda)=\phi(\lambda)$ for $\lambda \geq 0$. Moreover, since $\hphin(\lambda),\hphi(\lambda)>0$ for $\lambda>0$, one more application of Lemma \ref{lem:WH_factors_conv} gives, for every $\lambda \in (0,\delta)$,
    \[
    \lim_{n \to\infty}\phin(-\lambda) = \lim_{n \to \infty}\frac{\psin(\lambda)}{\hphin(\lambda)}= \frac{\psi(\lambda)}{\hphi(\lambda)} = \phi(-\lambda).
    \]
    Finally, we note with the help of Remark \ref{rem:H_levy_measure} that the integral representations of $\phin$ and $\phi$ are valid at least on $\C_{[-\delta,\infty])}$. By differentiating twice one concludes that we arrive at a sequence of concave smooth functions $\phin$ converging pointwise on $(-\delta,\infty)$ to a concave smooth function $\phi$. Then by \cite[Theorem 25.17]{Rockafellar97},
    \[
    \lim_{n \to \infty}(\phin)'(0) = \phi'(0) \quad \text{or equivalently} \quad \lim_{n \to \infty} \EE \Hn_1 = \EE H_1,
    \]
    as desired.
\end{proof}
\begin{remark}\label{rem:esscher2}
 Suppose $\EE e^{\beta \xi_1}<\infty$, $\limsup_{t \to \infty} \xi_t=\infty$ $\PP$-a.s., and let $\psie$ be the Laplace exponent of the L\'{e}vy process $\xie$ obtained from $\xi$ by the exponential change of measure (Esscher transform) \eqref{eq:cPP}, so that
 \[
    \psie(\lambda) = \psi(\lambda+\epsilon)-\psi(\epsilon), \quad z \in \C_{[-\epsilon,\beta-\epsilon]}.
 \]
 Assume $\psi(0)=0$, i.e. $\xi$ is not killed. Then according to Remark \ref{rem:Esscher_transform}, $\xie \to \xi$ as $\epsilon \to 0^+$ in the Skorokhod space. Moreover, we have $\lim_{t \to \infty} \xie_t=+\infty$ a.s. for $\epsilon \in (0,\beta)$ and one readily checks that $\psie(\lambda) \to \psi(\lambda)$ as $\epsilon \to 0^+$ for all $\lambda \in [0,\beta)$. Hence, this particular setting satisfies the assumptions of Proposition \ref{prop:pointwise} (and it will be used later in the proof of Theorem \ref{thm:moments}).
\end{remark}



\section{Proof of Theorem \ref{thm:moments}}\label{sec:moments_proof}
\subsection{Applications to Bernstein-Gamma functions}
Let us first specialise to the setting of \cite{PS18-ejp}.
Consider the Wiener-Hopf factorisation of the L\'{e}vy-Khintchine exponent given by
\begin{equation}\label{phipm}
\psi(z)-q = -\phi(-z)\hphi(z), \quad z \in i\R,
\end{equation}
where $q\geq 0$ is the killing rate of $\xi$, and $\phi,\hphi$ belong to the class of Bernstein functions; cf. \eqref{eq:WH}. Since the analytic properties of this class are at the heart of the approach in \cite{PS18-ejp}, we spend some time carefully introducing the notation and basic objects. The following excerpt is taken directly from \cite[Section 2 and Appendix A]{PS18-ejp}.

For any Bernstein function $\phi$, we let
\begin{align*}
\ga_\phi &= \inf \{u<0 \colon \phi \in {\rm A}_{(u,\infty)}\} \in [-\infty,0],\\
\gu_\phi &= \sup \{u \in [\ga_\phi,0]\colon \phi(u)=0\} \in [-\infty,0],\\
\oga_\phi &= \max \{\ga_\phi,\gu_\phi\} \in [-\infty,0],
\end{align*}
with the convention that $\sup \varnothing=-\infty$ and $\inf \varnothing=0$. For the sake of brevity, for the quantities above related to the Wiener-Hopf factors $\phi$, $\hphi$, we use the shortened notation $\ga,\gu,\oga$ and $\hga,\hgu,\hoga$, respectively. This convention is naturally transferred to other objects pertaining to $\phi$ and $\hphi$. We also (formally) define a generalised Weierstrass product as
\[
W_\phi(z) = \frac{e^{-\gamma_\phi z}}{\phi(z)} \prod_{k=1}^\infty \frac{\phi(k)}{\phi(z+k)}e^{\frac{\phi'(k)}{\phi(k)}z}, \quad z \in \C_+,
\]
where
\begin{equation*}
\gamma_\phi = \lim_{n \to \infty} \bigg( \sum_{k=1}^n \frac{\phi'(k)}{\phi(k)}- \ln \phi(n)\bigg) \in \bigg[ -\ln \phi(1)\,,\,\frac{\phi'(1)}{\phi(1)}-\ln \phi(1)  \bigg].
\end{equation*}
As it is shown in \cite{PS18-ejp, 2021-PP-MS-mams}, the function $W_\phi$ is 
the unique solution to the equation
\begin{equation}\label{eq:W_rec}
W_\phi(z+1) = \phi(z)W_\phi(z), \quad W_\phi(1)=1.
\end{equation}
As such, the (solutions of the) equations of the form above were studied by Webster in \cite{1997-Webster-jfaa}. Note that for the simple case $\phi(z)=z$ the function $W_\phi$ becomes the Gamma function. For this reason, $W_\phi$ is called in \cite{PS18-ejp} the Bernstein-Gamma function. According to \cite[Theorem 4.1]{PS18-ejp}, $W_\phi \in \rA_{(\oga_\phi,\infty)} \cap {\rm M}_{(\ga_\phi,\infty)}$ and $W_\phi$ is zero-free on $(\ga_\phi,\infty)$. Moreover, since $\phi$ is log-concave on $\R_+$, it follows from \cite{1997-Webster-jfaa} that $\R_+ \ni x \mapsto W_\phi(x)$ is log-convex on $\R_+$. For a future reference, we also let
    \[
    W_\phi(\ga_\phi):=\lim_{a \to (\ga_\phi)^+}\frac{W_\phi(1+a)}{\phi(a)},
    \]
    c.f. \eqref{eq:W_rec}, whenever the limit on the RHS is well defined. For an extensive analysis of the function $W_\phi$ and its properties, we refer the reader to \cite{PS18-ejp} and the references therein. 

Coming back to the setting of \eqref{phipm}, let us now assume that $\EE e^{\beta \xi_1}<\infty$ and $\phi(0)=0$, the latter implying $q=0$ and being equivalent to $\limsup_{t \to \infty} \xi_t=\infty$ $\PP$-a.s. Then by Proposition \ref{prop:kappa_extension} we get that $\phi \in \rA_{[-\beta,\infty)}$ and $\ga\leq-\beta$. By Remark \ref{rem:H_levy_measure} and differentiation under the integral one concludes that $\phi$ is strictly concave on the real half-line $[-\beta,\infty)$, and since $\phi(0)=0$, we get $\gu=0$ and consequently $\oga=0$.

Next, recall that $\psie(z) = \psi(z+\epsilon)-\psi(\epsilon)$ for $z \in \C_{[0,\beta-\epsilon]}$ and $\psie(0)=0$. We use the self-explanatory notation $\phie, \hphie$ and $\We,\hWe$ for objects pertaining to the L\'{e}vy process $\xi$ under $\PPa{\epsilon}$ with the characteristic exponent $\psie$. The classical Wiener-Hopf factorisation for $\xie$ reads
\[
    \psie(z) = -\phie(-z)\hphie(z) =-h_\epsilon(0)\kappae(0,-z)\hkappae(0-z), \quad z \in i\R.
\]
In view of Corollary \ref{cor:WH_changed} and \cite[Lemma A.1]{PS18-ejp}, we have $\ga_\epsilon \leq -\beta+\epsilon$ and one may take after recalling that $\ce \hce=h(\lambda_\epsilon)/h_\epsilon(0)$,
    \begin{align}\begin{aligned}\label{eq:cphi}
    \phie(z)&=h_\epsilon(0)\kappae(0,z)= h(\lambda_\epsilon)\kappa(\lambda_\epsilon,z-\epsilon), \quad &&z \in \C_{[-\beta+\epsilon,\infty)}, \\\hphie(z)&=\hkappae(0,z)=\hkappa(\lambda_\epsilon,z+\epsilon), \qquad &&z \in \C_{[0,\infty)},\end{aligned}
    \end{align}
    with $\psi(\lambda_\epsilon)=\epsilon$. Moreover, since $\lim_{t \to \infty} \xi_t=\infty$, $\PPa{\epsilon}$-a.s., we get as above that $\phie(0)=0$, so that $\gu_\epsilon=0$ and consequently, $\oga_\epsilon=0$, too.

    As we will shortly see, in our analysis it will be important to include the left endpoint of the domain of meromorphicity of $W$ in the application of results of \cite{2026-MM-MS-prob_surv}. Let us first collect some useful properties:
    \begin{lemma}\label{lem:phi_W_prop}
        Suppose that $\xi$ is not a compound Poisson process. Assume that $\EE e^{\beta \xi_1}<\infty$ and $\limsup_{t \to \infty}\xi_t=\infty$ $\PP$-a.s. Then the following hold:
        \begin{enumerate}[label={{\rm (\arabic*)}}]
        \item\label{lem:phi_W_prop_item_phipme2} $ \lim\limits_{\epsilon \to 0^+}\hphie(z) = \hphi(z)$ for all $ z \in \C_{[0,\infty)}$.
        \item\label{lem:phi_W_prop_item_phipme} $ \lim\limits_{\epsilon \to 0^+}\phie(z) = \phi(z)$ for all $z \in \C_{(-\beta,\infty)}$ and $\lim\limits_{\epsilon \to 0^+}\phie(z+ \epsilon) = \phi(z)$ for all $z \in \C_{[-\beta,\infty)}$.
        \item\label{lem:phi_W_prop_item_Wpme} For all $z \in \C_{(0,\infty)}$ we have $\lim\limits_{\epsilon \to 0^+} \We(z)=W(z)$ and $\lim\limits_{\epsilon \to 0^+} \hWe(z)=\hW(z)$.
        \item\label{lem:phi_W_prop_item_Wpme2} For all $x>0$ we have $\lim\limits_{\epsilon \to 0^+} \We(x+\epsilon)=W(x)$ and $\lim\limits_{\epsilon \to 0^+} \hWe(x-\epsilon)=\hW(x)$.
        \item\label{lem:phi_W_prop_item_We} $\lim_{\epsilon \to 0^+} \We(x+\epsilon) = W(x)$ for all $x \in [-\beta,\infty) \setminus \Z$.
        \end{enumerate}
    \end{lemma}
    \begin{proof}
    \ref{lem:phi_W_prop_item_phipme2}: according to \cite[Lemma A.1]{PS18-ejp} and \eqref{eq:cphi}, we have
    \[
    \hphi(z)=\hkappa(0,z) \quad \text{and} \quad \hphie(z)=\hkappa(\lambda_\epsilon,z+\epsilon) \quad \text{for} \quad z \in \C_{[0,\infty)}.
    \]
    Since $\lambda_\epsilon \to 0$ as $\epsilon \to 0^+$, the convergence $\hkappa(\lambda_\epsilon,z+\epsilon) \to \hkappa(0,z)$ as $\epsilon \to 0^+$ follows by the joint continuity of $\hkappa$.
    
    \ref{lem:phi_W_prop_item_phipme}:  \cite[Lemma A.1]{PS18-ejp}, Proposition \ref{prop:kappa_extension} and \eqref{eq:cphi} imply
    \begin{equation*}
    \phi(z)=h(0)\kappa(0,z) \quad z \in \C_{[-\beta,\infty)}
    \end{equation*}
    and
    \begin{equation*}
    \phie(z) = h(\lambda_\epsilon)\kappa(\lambda_\epsilon,z-\epsilon) \quad \text{for} \quad z \in \C_{[-\beta+\epsilon,\infty)}.
    \end{equation*}
    By the monotone convergence theorem, $h(\lambda_\epsilon) \to 1=h(0)$ as $\epsilon \to 0^+$. Let us consider three cases. If $\Re z >0$, then, for $\epsilon <\Re z$ we have $\kappa(\lambda_\epsilon,z-\epsilon) \to \kappa(0,z)$ as $\epsilon \to 0^+$ by the joint continuity of $\kappa$. For $z=ix$ with $x \in \R$ we recover L\'{e}vy-Khinchine exponents of $\He$ and $H$, respectively, so the convergence follows by Remark \ref{rem:Esscher_transform} and Lemma \ref{lem:WH_factors_conv}. For $\Re z \in \C_{(-\beta,0)}$ we use the extended Wiener-Hopf factorisation together with Remark \ref{rem:esscher2} and \ref{lem:phi_W_prop_item_phipme2} to get the first statement. Since $\phie(z+\epsilon)=h(\lambda_\epsilon)\kappa(\lambda_\epsilon,z)$, the second part follows by an even simpler argument.
    \bigskip

    \ref{lem:phi_W_prop_item_Wpme}: follows readily from \cite[Lemma 4.10]{PS18-ejp}, with the pointwise limits secured by \ref{lem:phi_W_prop_item_phipme2} and \ref{lem:phi_W_prop_item_phipme}.

    \ref{lem:phi_W_prop_item_Wpme2}: follows from the fact that the pointwise limit of log-convex, hence convex functions is uniform on every compact subinterval, see e.g. Rockafellar \cite[Theorem 10.8]{Rockafellar97}.

    \ref{lem:phi_W_prop_item_We}: for $x \in [-\beta,0) \setminus \Z$ we note that $\phie(0)=0$, so by concavity, $\phie(x)<0$ for $x \in [\ga_\epsilon,0)$. Thus, letting $\epsilon>0$ small enough so that $x+\epsilon\notin\Z$, we can use the recurrent equation to get
    \[
    \We(x+\epsilon) = \frac{\We(x-[x]+\epsilon)}{\phie(x+\epsilon)\phie(x+1+\epsilon)\ldots\phie(x-[x]+\epsilon)}.
    \]
    The statement now follows from \ref{lem:phi_W_prop_item_phipme} and \ref{lem:phi_W_prop_item_Wpme2}.
\end{proof}

\subsection{Proof}
For the proof of Theorem \ref{thm:moments}, a careful analysis of the infinitesimal generator of the underlying L\'{e}vy process and the application of the Dynkin formula is crucial. Since its direct application to the exponential function is not possible due to the fact that this function is not in the domain of this generator, instead we apply the approximation argument using appropriately chosen sequence of smooth, compactly supported functions. To this end, we let  $\chi_1 \colon \R \mapsto \R$ be an auxiliary $\Cs(\R)$ function satisfying $\chi_1(x)=0$ for $x\leq 0$ and $\chi_1(x)=1$ for $x\geq 1$. Then we set $f(x)=\chi_1(x)e^{\beta x}$, so that
	\begin{enumerate}[label={\arabic*.}]
        \item $f \in \Cs(\R)$.
		\item $f(x)=0$ for $x \leq 0$ and $f(x)=e^{\beta x}$ for $x \geq 1$.
        \item $f(x), |f'(x)|,|f''(x)| \leq ce^{\beta x}\ind_{\{x>0\}}$.
	\end{enumerate}

Let us start with an elementary observation, which is presented in the form of the lemma for convenience of the future reference later on. We stress that the constant $c$ is independent of $n$, and thus the established bounds are uniform and may (and will) serve as a majorising functions in the dominated convergence theorem.
\begin{lemma}\label{lem:fn_seq}
	There exists a sequence $(f_n)_{n \in \N}$ of real functions and a constant $c>0$ depending only on $\beta$ such that:
	\begin{enumerate}[label={\rm (\arabic*)}]
		\item $f_n \in \Csc(\R)$ with $f_n(x)=f(x)$ for $x \leq n$ and $f_n(x) = 0$ for $x \geq n+1$, and for all $x \in \R$ we have $f_n(x) \nearrow f(x)$, $f_n'(x) \to f'(x)$ and $f_n''(x) \to f''(x)$ as $n \to \infty$.
		\item For all $x \in \R$ we have $f_n(x),|f_n'(x)|, |f_n''(x)|\leq ce^{\beta x}\ind_{\{x>0\}}$.
	\end{enumerate}
\end{lemma}
\begin{proof}
    Let $\chi_2 \colon \R \mapsto \R$ be a $\Cs(\R)$ satisfying $\chi_2(x)=1$ for $x \leq 0$ and $\chi_2(x)=0$ for $x \geq 1$. We then let $f_n(x)=\chi_2(x-n)f(x)$, so that $f_n(x)=f(x)$ for $x \leq n$ and we easily see that all the conditions of the lemma are satisfied.
\end{proof}
Recall that for $\lambda \geq 0$ and $\beta \in \R$ we denote
\[
\Vo{\beta}(\lambda) = \EE \bigg[ \int_0^{\zeta} e^{-\lambda s} e^{\beta \xi_{\tau_s}}\ud s\bigg], \quad \lambda \geq 0,\,\beta \in \R.
\]
Let $\calV = \{\beta \in \R \colon \Vo{\beta}(\lambda)<\infty \text{ for } \lambda>0\}$.
Note that $\calV$ is a convex set (i.e. an interval, possibly a singleton) that contains the origin.
Indeed, if $\beta_1,\beta_2 \in \calV$ and $\beta_1<\beta<\beta_2$, then, for $\lambda>0$,
\begin{align*}
\Vo{\beta}(\lambda)&=\EE \bigg[ \int_0^{\zeta} e^{-\lambda s} e^{\beta\xi_{\tau_s}} \ind_{\{\xi_{\tau_s} < 0\}}\ud s \bigg] + \EE \bigg[ \int_0^{\zeta} e^{-\lambda s} e^{\beta\xi_{\tau_s}} \ind_{\{\xi_{\tau_s} \geq 0\}}\ud s \bigg] \\ &\leq
\EE \bigg[ \int_0^{\zeta} e^{-\lambda s} e^{\beta_1\xi_{\tau_s}} \ind_{\{\xi_{\tau_s} < 0\}}\ud s \bigg] + \EE \bigg[ \int_0^{\zeta} e^{-\lambda s} e^{\beta_2\xi_{\tau_s}} \ind_{\{\xi_{\tau_s}\geq 0\}}\ud s \bigg]  \\ &\leq V_{\beta_1}(\lambda)+\Vo{\beta_2}(\lambda),
\end{align*}
so that $\beta \in \calV$ as claimed.

In what follows, we denote by $\check\xi$ a non-killed version of $\xi$ obtained by allowing $\xi$ to evolve for infinite time horizon. Accordingly, $\cA$ and $\cpsi$ are the exponential functional and Laplace exponent associated with $\cxi$, respectively. Thus, $\cpsi(z)=\psi(z)+q$ for $z \in \rAC$, where $q$ is the killing intensity of $\xi$.

\begin{lemma}\label{lem:Vq}
	If for some $\beta\in\R$ and $\lambda>0$, $e^{\psi(\beta)}=\EE[e^{\beta\xi_1}]<\infty$ and $\Vo{\beta-1}(\lambda)<\infty$,
	then $\Vo{\beta}(\lambda)<\infty$ and
	\begin{equation}\label{eq:23}
		\lambda \Vo{\beta}(\lambda) \leq \cpsi(\beta)\Vo{\beta-1}(\lambda)+1.
	\end{equation}
	Moreover, if $ e^{\psi(\beta+1)}=\EE e^{(\beta+1)\xi_1}<\infty$, then we have the equality in \eqref{eq:23}.
\end{lemma}
Note that $\Vo{\beta}(\lambda)$ and $\Vo{\beta-1}(\lambda)$ in the statement above correspond to the possibly killed $\xi$, while $\cpsi$ is associated with its conservative version $\cxi$.
\begin{proof}
	Assume first that $\EE e^{\beta \xi_1}<\infty$. By change of variables and writing $\cxi$ instead of $\xi$,
	\begin{align}
		\Vo{\beta}(\lambda) &= \EE \bigg[\int_0^\infty e^{-qt}  e^{-\lambda \cA_t}e^{(1+\beta)\cxi_t}\ud t\bigg]\nonumber \\ &= -\frac{1}{\lambda}\EE \bigg[\int_0^\infty e^{-qt}\bigg(\frac{\od}{\od t}e^{-\lambda \cA_t}\bigg) e^{\beta\cxi_t}\ud t\bigg]\nonumber \\ &=
\frac{1}{\lambda} \EE \bigg[\int_0^\infty e^{-qt} \lim_{h \to 0^+} \frac{e^{-\lambda \cA_{t-h}}-e^{-\lambda \cA_t}}{h}\ind_{\{h<t\}}e^{\beta\cxi_t}\ud t\bigg]\nonumber \\ &\leq \frac{1}{\lambda} \liminf_{h \to 0^+} \EE\bigg[ \int_h^\infty e^{-qt} \frac{e^{-\lambda \cA_{t-h}}-e^{-\lambda \cA_t}}{h}e^{\beta \cxi_t}\ud t\bigg],\label{eq:12}
	\end{align}
	where in the last inequality we use the Fatou lemma and the fact that the function $t \mapsto A_t$ is non-decreasing.
	
	Let us analyse the (non-negative, but possibly infinite) right-hand side of \eqref{eq:12}. Recall that
	\[
	\Vo{\beta-1}(\lambda) = \EE\bigg[\int_0^\infty e^{-qt}e^{-\lambda \cA_t}e^{\beta\cxi_t}\ud t\bigg] < \infty.
	\]
	Thus, 
	\begin{align}
		\EE \bigg[\int_h^\infty e^{-qt} \frac{e^{-\lambda \cA_{t-h}}-e^{-\lambda \cA_t}}{h}e^{\beta\cxi_t}\ud t\bigg] &= \EE \bigg[\frac1h\bigg(\int_h^\infty e^{-qt} e^{-\lambda \cA_{t-h}}e^{\beta\cxi_t}\ud t-\int_h^\infty e^{-qt} e^{-\lambda \cA_t}e^{\beta\cxi_t}\ud t\bigg)\bigg] \nonumber \\ &=
 \EE \bigg[\frac1h\bigg(\int_h^\infty e^{-qt}e^{-\lambda \cA_{t-h}}e^{\beta(\cxi_t-\cxi_{t-h}+\cxi_{t-h})}\ud t-\int_h^\infty
 e^{-qt} e^{-\lambda \cA_t}e^{\beta\cxi_t}\ud t\bigg)\bigg] \nonumber \\ &= \frac1h\bigg( \Vo{\beta-1}(\lambda) \Big(e^{h \cpsi(\beta)}-1 \Big)+\EE \bigg[\int_0^h e^{-qt} e^{-\lambda \cA_t}e^{\beta\cxi_t}\ud t\bigg]\bigg), \label{eq:22}
	\end{align}
	where the last equality is justified by the following fact. For $h<1$ we clearly have $\int_0^h e^{-qt} e^{-\lambda \cA_t}e^{\beta\cxi_t}\od t \leq he^{\beta\overline{\cxi}_1}$ and our assumption implies that $\EE e^{\beta\overline{\cxi}_1}<\infty$. Indeed, if $\e$ is an independent exponentially distributed time with parameter, say, $\lambda>0$, then $\overline{\cxi}_\e \geq 0$ and $\cxi_{\e}-\overline{\cxi}_\e \leq 0$ are independent, so that
$\EE e^{\beta \cxi_{\e}} = \EE \big[ e^{\beta(\cxi_{\e}-\overline{\cxi}_{\e})}\big] \EE \big[ e^{\beta\overline{\cxi}_{\e}} \big]<\infty$ for $\lambda >\psi(\beta)$; see e.g. \cite[Theorem 5, p. 159]{Bertoin96}. Moreover,
$\infty>\EE \big[ e^{\beta\overline{\cxi}_\e} \big]\geq \lambda \int_1^\infty e^{-\lambda s} \EE \big[ e^{\beta\overline{\cxi}_{s}} \big]\ud s \geq e^{-\lambda} \EE \big[ e^{\beta\overline{\cxi}_{1}} \big]$, which completes the proof. The above also justifies the application of the dominated convergence theorem, so that
	\[
	\lim_{h \to 0^+} \Bigg(\frac{1}{h} \EE \bigg[\int_0^h e^{-qt} e^{-\lambda \cA_t}e^{\beta\cxi_t}\od t\bigg] \Bigg) = 1.
	\]
	It follows that the right-hand side of \eqref{eq:22}, and consequently the right-hand side of \eqref{eq:12}, has a limit as $h \to 0^+$ with
	\[
	\lim_{h \to 0^+} \EE\bigg[ \int_h^\infty e^{-qt}\frac{e^{-\lambda \cA_{t-h}}-e^{-\lambda \cA_t}}{h}e^{\beta\cxi_t}\ud t\bigg] =  \cpsi(\beta)\Vo{\beta-1}(\lambda)+1,
	\]
	and \eqref{eq:23} follows.
	
	Now assume that $\EE e^{(\beta+1)\xi_1}<\infty$. After change of variables $t-h=u$, by the mean-value theorem, for almost every $(\omega,t) \in \Omega\times \Rp$ we may bound the integrand on the right-hand side of \eqref{eq:12} as follows:
	\begin{align*}
	\frac{1}{\lambda} \frac{e^{-\lambda \cA_{t}}-e^{-\lambda \cA_{t+h}}}{h}e^{\beta\cxi_{t+h}} &\leq \frac{1}{\lambda}\sup_{s \in [t,t+h]} \bigg|\frac{\od}{\od s}e^{-\lambda \cA_s}\bigg|e^{\beta\cxi_{t+h}} \\ &= \sup_{s \in [t,t+h]} \bigg(e^{-\lambda \cA_s}e^{\cxi_s}\bigg)e^{\beta\cxi_{t+h}} \\ &\leq e^{-\lambda \cA_t}e^{(\beta+1)\cxi_t} \cdot e^{\beta(\cxi_{t+h}-\cxi_t)} \sup_{s \in [t,t+h]} e^{\cxi_s-\cxi_t} \\ &\leq e^{-\lambda \cA_t}e^{(\beta+1)\cxi_t} \cdot e^{(\beta+1)\overline{\cxi}_1^t},
	\end{align*}
	where $\overline{\cxi}_1^t = \sup_{s \in [0,1]}(\cxi_{t+s}-\cxi_t)$. Since $\EE e^{(\beta+1)\cxi_1}<\infty$, by the same argument as above we get that $\EE e^{(\beta+1)\overline{\cxi}_1}<\infty$. Thus, by the independence of increments of $\cxi$,
the fact that $\overline{\cxi}_1^t \ed \overline{\cxi}_1$ for every $t \geq 0$, and the first part of the proof, we conclude that
	\[
	\EE \bigg[ \int_0^\infty e^{-qt}  e^{-\lambda \cA_t}e^{(\beta+1)\cxi_t} e^{(\beta+1)\overline{\cxi}_1^t} \ud t\bigg] \leq \Vo{\beta}(\lambda) \EE e^{(\beta+1)\overline{\cxi}_1}<\infty,
	\]
	and one may apply the dominated convergence theorem to get the equality in \eqref{eq:12}.
\end{proof}

\begin{proposition}\label{prop:V_q_negative}
	Suppose that either $\lim_{t \to \infty }\xi_t=-\infty$ $\PP$-a.s. or $\xi$ is killed with rate $q>0$. If $\beta > 0$ and $e^{\psi(\beta)}=\EE e^{\beta\xi_1}<\infty$, then $\Vo{\beta-1}(\lambda)<\infty$ for all $\lambda>0$ and consequently, $\Vo{\beta}(\lambda)<\infty$ for all $\lambda>0$.
\end{proposition}
\begin{proof}
    Assume first that $\beta \in (0,1]$. By assumption, $\psi(s)$ is well defined, at least for all $s \in [0,\beta]$. If $\lim_{t \to \infty}\xi_t=-\infty$ a.s., then we have that $(\psi')(0^+) < 0$ and consequently, $\psi$ is strictly negative around the origin. If $\xi$ is killed with $q>0$, then $\psi(0)=-q<0$ and one arrives at the same conclusion. We consider two cases. If $\psi(\beta)<0$, then, recalling from the proof of Lemma \ref{lem:Vq} that $\psi(\beta)=\check{\psi}(\beta)-q$, we simply write
	\begin{align*}
		\Vo{\beta-1}(\lambda) = \EE \bigg[ \int_0^\infty e^{-qt} e^{-\lambda \check{A}_t}e^{\beta \check{\xi}_t}\ud t \bigg] \leq \int_0^\infty e^{(\check{\psi}(\beta)-q)t}\od t = \frac{1}{-\psi(\beta)}<\infty,
	\end{align*}
	and $\Vo{\beta}(\lambda)<\infty$ follows by Lemma \ref{lem:Vq}. If $\psi(\beta) \geq 0$, then there is $s_0 \in (0,\beta)$ such that $\psi(s_0)<0$, and by the first part of the proof $\Vo{s_0-1}(\lambda)<\infty$.
Lemma \ref{lem:Vq} entails $\Vo{s_0}(\lambda)<\infty$, so by convexity, $\Vo{\beta-1}(\lambda)<\infty$, and one more application of Lemma \ref{lem:Vq} yields the claim. If $\beta  \in (1,2]$, then clearly $\EE e^{(\beta-1)\xi_1}=e^{\psi(\beta-1)}<\infty$ and $\Vo{\beta-1}(\lambda)<\infty$ by the previous line.
Lemma \ref{lem:Vq} then implies that $\Vo{\beta}(\lambda)<\infty$ as well. Repeating this recursive procedure $\lceil \beta\rceil-1$ times we get that $\Vo{\beta}(\lambda)<\infty$ for any $\beta>0$, provided that $\EE e^{\beta\xi_1}=e^{\psi(\beta)}<\infty$.
\end{proof}
\begin{remark}\label{remarkproposition}
    If $\lim_{t \to \infty}\xi_t=\infty$ $\PP$-a.s, then it holds that $\Vo{-1}(\lambda)<\infty$ for all $\lambda>0$. Indeed, we note that
    \begin{align*}
		\Vo{-1}(\lambda) = \EE \bigg[ \int_0^\infty e^{-\lambda A_t} \ud t \bigg].
	\end{align*}
	By choosing any $a \in (0,1)$ and using \cite[Theorem 3.9]{2026-MM-MS-transAMS} with $F(x)=e^{-\lambda x}$ (note that $I_\xi(t)=\int_0^te^{-\xi_s}\ud s$ in the notation of \cite{2026-MM-MS-transAMS}, see the Introduction of \cite{2026-MM-MS-transAMS}, so the condition $\lim_{t \to \infty} \xi_t = -\infty$ a.s. is satisfied) we get the conclusion. Since $\Vo{0}(\lambda)=\frac{1}{\lambda}$, by convexity we get that $\Vo{\beta}(\lambda)<\infty$ for all $\beta \in [-1,0]$. It follows readily that if $\EE e^{\beta \xi_1}<\infty$, then by applying analogous recursive procedure as in the proof of Proposition \ref{prop:V_q_negative} one may conclude that $\Vo{\beta}(\lambda)<\infty$ for all $\lambda > 0$.
\end{remark}

Recall from \cite[Corollary 2.9]{LevyMatters3} that the infinitesimal generator $\genxi$ of $\xi$ has the following form:
\begin{align}\label{eq:genX}
	\genxi g(x) &=-q g(x)+ bg'(x)+\frac{\sigma^2}{2}g''(x) + \int_{\R} \big(g(x+y)-g(x)-yg'(x)\ind_{\{|y| \leq 1\}}\big) \,\nu(\od y),
\end{align}
for appropriate, say $\Csc(\R)$ function $g$.
Now we are ready to prove Theorem \ref{thm:moments}.
\begin{proof}[Proof of Theorem \ref{thm:moments}]
	Let us start with a useful observation. Fix arbitrary $t >0$ and let $(f_n)_{n \in \N}$ be the sequence of functions approximating $f$ constructed in Lemma \ref{lem:fn_seq}. Since $f_n \in \Csc(\R)$ for each $n \in \N$, it is in the domain of the generator $\genxi$ of $\xi_t$ by virtue of \cite[Corollary 2.9]{LevyMatters3}. It follows then (see e.g. \cite[Theorem 1.36]{LevyMatters3}) that the process
    \[
    M_u := f_n(\xi_u)-f_n(\xi_0)-\int_0^u \genxi f_n(\xi_s)\ud s
    \]
    is a martingale with respect to the natural filtration of $\xi$. Let $m \in \N$. Proceeding as in \cite[proof of Corollary 4.2]{LevyMatters3} we find that
    \begin{align*}
        f_n(\xi_{\tau_t \wedge m}) - f_n(\xi_0)-M_{\tau_t \wedge m} &= \int_0^{\tau_t \wedge m} \genxi f_n(\xi_s)\ud s \\ &= \int_0^{\tau_t \wedge m} e^{-\xi_s} \genxi f_n(\xi_s) e^{\xi_s} \ud s \\ &= \int_0^{t  \wedge A_m \wedge \zeta} e^{-\xi_{\tau_s}} \genxi f_n(\xi_{\tau_s}) \ud s.
    \end{align*}
    Since $\tau_t \wedge m$ is a bounded stopping time, we have from the optional stopping theorem that
    \begin{equation}\label{eq:gen_eq}
    \EE f_n(\xi_{\tau_t \wedge m}) = \EE \bigg[\int_0^{t \wedge A_m \wedge \zeta} e^{-\xi_{\tau_s}} \genxi f_n(\xi_{\tau_s})\ud s \bigg].
    \end{equation}
    We split the expectation on the RHS of \eqref{eq:gen_eq} as follows:
    \begin{align*}
        \EE \bigg[\int_0^{t \wedge A_m \wedge \zeta} e^{-\xi_{\tau_s}} \genxi f_n(\xi_{\tau_s})\ud s \bigg] = \rI_1(t,n,m) + \rI_2(t,n,m) - \rI_3(t,n,m) + \rI_4(t,n,m) + \rI_5(t,n,m),
    \end{align*}
    where, according to \eqref{eq:genX},
    \begin{align*}
        \rI_1(t,n,m) &=  \EE \bigg[\int_0^{t \wedge A_m \wedge \zeta} e^{-\xi_{\tau_s}} \bigg(-qf_n(\xi_{\tau_s}) + bf'_n(\xi_{\tau_s}) + \frac{\sigma^2}{2}f''_n(\xi_{\tau_s})\bigg) \ud s\\
        \rI_2(t,n,m) &= \EE \bigg[\int_0^{t \wedge A_m \wedge \zeta} e^{-\xi_{\tau_s}} \int_{|y|\leq 1} \big(f_n(\xi_{\tau_s}+y)-f_n(\xi_{\tau_s})-yf'_n(\xi_{\tau_s})\big)\,\nu(\od y) \ud s\bigg],\\
        \rI_3(t,n,m) &= \EE \bigg[\int_0^{t \wedge A_m \wedge \zeta} e^{-\xi_{\tau_s}} f_n(\xi_{\tau_s}) \int_{|y|>1}\nu(\od y) \ud s\bigg],\\
        \rI_4(t,n,m) &= \EE \bigg[\int_0^{t \wedge A_m \wedge \zeta} e^{-\xi_{\tau_s}} \int_{(-\infty,-1)} f_n(\xi_{\tau_s}+y) \,\nu(\od y) \ud s\bigg],\\ \intertext{and}
        \rI_5(t,n,m) &= \EE \bigg[\int_0^{t \wedge A_m \wedge \zeta} e^{-\xi_{\tau_s}} \int_{(1,\infty)} f_n(\xi_{\tau_s}+y) \,\nu(\od y) \ud s\bigg].
    \end{align*}
    Assume first that $\beta\leq 1$. We claim that the integrals $\rI_1$ -- $\rI_4$ (are absolutely convergent and) converge to a finite limit as $m,n \to \infty$, irrespectively of the order of taking limits. Indeed, $\rI_1$ follows from Lemma \ref{lem:fn_seq} together with the facts that $\beta \in (0,1]$ and $\supp f_n \subseteq [0,\infty)$, so that
    \begin{align}\label{eq:14}\begin{aligned}
    \ind_{(0,t \wedge A_m \wedge \zeta})(s)\Bigg| e^{-\xi_{\tau_s}} \bigg(- q f_n(\xi_{\tau_s}) + bf'_n(\xi_{\tau_s}) + \frac{\sigma^2}{2}f''_n(\xi_{\tau_s})\bigg)\Bigg| &\leq c\ind_{\{s<t \wedge \zeta}\} e^{-(1-\beta)\xi_{\tau_s}} \ind_{\{\xi_{\tau_s}>0\}} \\ &\leq c\ind_{\{s<t\}}.\end{aligned}
    \end{align}
    Using the fact that $\nu((-1,1)^c)<\infty$, the same argument applies to $\rI_3$. For $\rI_2$ we use the standard bound
    \begin{align*}
    &\ind_{\{|y|\leq 1\}}\ind_{(0,t \wedge A_m \wedge \zeta})(s)e^{-\xi_{\tau_s}}\Big|f_n(\xi_{\tau_s}+y)-f_n(\xi_{\tau_s})-yf'_n(\xi_{\tau_s})\Big| 
    \\ &\qquad \leq e^{-\xi_{\tau_s}}|y|^2\ind_{\{|y|\leq 1\}}
    \sup_{\eta \in B(\xi_{\tau_s},1)} f_n''(\eta) \ind_{\{s<t \wedge \zeta\}} \\ &\qquad\leq c|y|^2\ind_{\{|y|\leq 1\}}e^{-(1-\beta)\xi_{\tau_s}}\ind_{\{\xi_{\tau_s}>-1\}}\ind_{\{s<t \wedge \zeta\}} \\ &\qquad \leq c|y|^2\ind_{\{|y|\leq 1\}}\ind_{\{s<t\}}.
    \end{align*}
    Finally, for $\rI_4$ we have
    \begin{align}\label{eq:16}\begin{aligned}
    \ind_{(0,t \wedge A_m \wedge \zeta})(s)\ind_{\{y<-1\}}e^{-\xi_{\tau_s}}  f_n(\xi_{\tau_s}+y) &\leq e^{\beta y}\ind_{\{y<-1\}}e^{-(1-\beta)\xi_{\tau_s}}\ind_{\{\xi_{\tau_s}>0\}}\ind_{\{s<t \wedge \zeta\}}  \\ &\leq c\ind_{\{y<-1\}}\ind_{\{s<t\}},\end{aligned}
    \end{align}
    and the claim is proved by the dominated convergence. Next, using first the quasi-left continuity of $\xi$ and continuity and boundedness of $f_n$ and then the monotone convergence theorem, we get that
    \[
    \lim_{n \to \infty} \lim_{m \to \infty} \EE f_n(\xi_{\tau_t \wedge m}) = \EE f(\xi_{\tau_t}),
    \]
    although the limit may be infinite. Since the limit of $\rI_5$ also exists by the monotone convergence irrespectively of the order of taking limits, we may take the limit in \eqref{eq:gen_eq} first with $m \to \infty$ and then with $n \to \infty$ to conclude that
    \begin{equation}\label{eq:moments_id}
    \EE f(\xi_{\tau_t}) < \infty \qquad \text{if and only if} \qquad \EE \bigg[ \int_0^{t \wedge \zeta}e^{-\xi_{\tau_s}}\int_{(1,\infty)} f(\xi_{\tau_s}+y)\,\nu(\od y)\ud s \bigg]< \infty.
    \end{equation}
    Let us stress that the equivalence above holds for all L\'{e}vy processes $\xi$, irrespectively of the presence of killing or their behaviour at infinity.

    Now we proceed to the main part of the proof. Suppose first that $\limsup_{t \to \infty} \xi_t=\infty$ $\PP$-a.s. and that $\xi$ is not a compound Poisson process. We will first prove that $\ref{thm_moments_item_levy} \Rightarrow \ref{thm_moments_item_entrance} \Rightarrow \ref{thm_moments_item_V0} \Rightarrow \ref{thm_moments_item_Vx} \Rightarrow \ref{thm_moments_item_ssmpp} \Rightarrow \ref{thm_moments_item_levy}$, where
    \begin{enumerate}[label={{\rm (\arabic*')}}]
    \item\label{thm_moments_item_ssmpp} For all $x>0$, we have $\mathbb{E}_x X _t^\beta<\infty$ for almost all $t>0$.
    \end{enumerate}
    Recall that in this case we have $\zeta=\infty$, $\P_x$-a.s., for all $x>0$. \medskip
    
    $\ref{thm_moments_item_levy} \Rightarrow \ref{thm_moments_item_entrance}$: Note that from the independence and homogeneity of the increments of $\xi$ and the scaling property under $\mathbb{P}_0$, we can assume that $t=1$. We know from Proposition \ref{prop:kappa_extension} and Remark \ref{rem:H_levy_measure} that $\EE e^{\beta H_1}<\infty$, where $H$ is the upward ladder height process associated with $\xi$, and in particular $\EE H_1<\infty$. By \cite[Theorem 4]{2012-LC-AEK-JCP-VR-aop}, the pssMp $(X,\P_x)$ admits the entrance law at $0$ and the measure $\P_0$ is well defined.

    Let $\xie$ be a L\'{e}vy process with the Laplace exponent $\psie(\lambda) = \psi(\lambda+\epsilon)-\psi(\epsilon)$ for $\Re \lambda \in [0,\beta-\epsilon]$, and denote by $\Xe$ the pssMp associated with $\xie$. Proposition \ref{prop:pointwise} asserts that the assumption \ref{A:EH_conv} is satisfied for $\He$. Thus, by Theorem \ref{thm:ssMps_0_conv} and the Fatou lemma we get
    \[
    \E_0 X_1^\beta \leq \liminf_{\epsilon \to 0^+} \E_0 \bigg[\Big(\Xe_1\Big)^{\beta-\epsilon}\bigg].
    \]
    Assume first that $\beta \in (0,1)$ and recall that $\gu_\epsilon=0$ and $\ga_\epsilon \leq -\beta+\epsilon$ in the notation of \cite{PS18-ejp}, see the discussion after \eqref{eq:cphi}. With the aim to use \cite[Theorem 3.1 and Corollary 3.4]{2026-MM-MS-prob_surv} in the left endpoint, we observe from the recurrent equation \eqref{eq:W_rec} that the expression
    \begin{equation}\label{eq:Wp_limit}
    \lim_{a \to (-\beta+\epsilon)^+}\frac{\Gamma(a)}{\We(a)}  =\lim_{a \to (-\beta+\epsilon)^+} \frac{\phie(a)}{a} \frac{\Gamma(1+a)}{\We(1+a)} = \frac{\phie(-\beta+\epsilon)}{-\beta+\epsilon} \frac{\Gamma(1-\beta+\epsilon)}{\We(1-\beta+\epsilon)}
    \end{equation}
    is well defined, c.f. \cite[p. 138]{2026-MM-MS-prob_surv}. Thus, by \cite[Theorem 1]{02-JB-MY-potanal} and \cite[Theorem 3.1 and Corollary 3.4]{2026-MM-MS-prob_surv},
    \[
    \E_0\Big[ \Big(\Xe_1\Big)^{\beta-\epsilon}\Big] = \frac{1}{m_\epsilon}\EEa{\epsilon} \bigg[ \frac{1}{\rI^{1+\beta-\epsilon}} \bigg] = \frac{1}{m_\epsilon} \hphie(0) \frac{\phie(-\beta+\epsilon)}{-\beta+\epsilon} \frac{\Gamma(1-\beta+\epsilon)}{\We(1-\beta+\epsilon)} \hWe(1+\beta-\epsilon),
    \]
    where $I=\int_0^\infty e^{\xie_t}\ud t$.
    Recall that the (extended) Wiener-Hopf factorisation yields $\psie(\lambda)=-\phie(-\lambda)\hphie(\lambda)$ for $\Re \lambda \in [0,\beta-\epsilon]$.
    Since $\xie$ is not killed and drifts to $\infty$, we have $\psie(0) = \phie(0)=0$, so dividing this equality by $\lambda$ and taking the limit as $\lambda \to 0^+$ yields
    \[
    (\psie')(0^+) = m_\epsilon = (\phie)'(0)\hphie(0) = \hphie(0)\EE \He_1.
    \]
    Recall that \ref{A:EH_conv} holds in this setting, so it follows that
    \[
    \lim_{\epsilon \to 0^+} \frac{\hphie(0)}{m_\epsilon} = \lim_{\epsilon \to 0^+} \frac{1}{\EE \He_1}=\frac{1}{\EE H_1}<\infty.
    \]
    The application of Lemma \ref{lem:phi_W_prop} then yields
    \[
    \liminf_{\epsilon \to 0^+} \E_0\bigg[ \Big(\Xe_1\Big)^{\beta-\epsilon}\bigg] \leq \frac{1}{\EE H_1} \frac{\phi(-\beta)}{-\beta} \frac{\Gamma(1-\beta)}{W(1-\beta)}\hW(1+\beta)<\infty,
    \]
    and \ref{thm_moments_item_entrance} follows in this case.

    Now let $\beta=1$. Then the identity \eqref{eq:Wp_limit} can be further rewritten as
    \[
    \lim_{a \to (-1+\epsilon)^+}\frac{\Gamma(a)}{\We(a)} = \frac{\phie(-1+\epsilon)}{-1+\epsilon} \frac{\Gamma(\epsilon)}{\We(\epsilon)} = \frac{\phie(-1+\epsilon)}{-1+\epsilon} \frac{\phie(\epsilon)}{\epsilon} \frac{\Gamma(1+\epsilon)}{\We(1+\epsilon)}.
    \]
    Note that by the concavity of $\phie$,
    \[
    \limsup_{\epsilon \to 0^+}\frac{\phie(\epsilon)}{\epsilon} \leq \limsup_{\epsilon \to 0^+}(\phie)'(0)=\limsup_{\epsilon \to 0^+} \EE \He_1=\EE H_1<\infty,
    \]
    and \ref{thm_moments_item_entrance} follows by the same argument as above.

    Finally, if $\beta > 1$, then from the proof of \cite[Proposition 2]{02-JB-MY-toulouse} we get that
    \[
    \E_0 X_1^\beta = \frac{\psi(\beta)}{\beta} \E_0 X_1^{\beta-1},
    \]
    and \ref{thm_moments_item_entrance} follows from the recursive argument and the already proved case $\beta \in (0,1]$.

    $\ref{thm_moments_item_entrance} \Rightarrow \ref{thm_moments_item_V0}$: by the scaling property, 
    $\E_0 X_t^\beta=ct^\beta$ with $c=\E_0 X_1^\beta \in (0,\infty)$. Thus, for every $\lambda>0$,
    \[
    \Vx{\beta}{0}(\lambda) = \E_0\bigg[ \int_0^{\infty} e^{-\lambda t}X_t^\beta \ud t\bigg] = c\int_0^\infty e^{-\lambda t}t^\beta \ud t =c\lambda^{-(1+\beta)}\Gamma(1+\beta)< \infty.
    \]

    $\ref{thm_moments_item_V0} \Rightarrow \ref{thm_moments_item_Vx}$: by the Markov property, for every $\lambda>0$ and $t > 0$,
    \begin{align}\label{eq:V0-MP}
        \infty>\Vx{\beta}{0}(\lambda) \geq \E_0 \bigg[\int_t^{\infty} e^{-\lambda s}X_s^\beta\ud s \bigg]= e^{-\lambda t}\int_0^{\infty} e^{-\lambda s}\E_{X_t}\big[X_s^\beta\big] \ud s=e^{-\lambda t} \E_0 \Big[\Vx{\beta}{X_t}(\lambda)\Big].
    \end{align}
    Let us note that from the scaling property, for all $0<a<b<\infty$ there is $t>0$ such that $\P_0(X_t \in (a,b))>0$. Indeed, if there were $a,b$ such that for all $t>0$,
    \[
    \P_0\big(X_t \in (a,b)\big) = \P_0 \big( X_1 \in (t^{-1}a, t^{-1}b)\big)=0,
    \]
    then by choosing properly a countable number  of values of $t$, one would obtain $\P_0\big(X_1 \in (0,\infty)\big)=0$, which is impossible since $X_1$ is non-degenerate under $\P_0$.

    Thus, we get from \eqref{eq:V0-MP} that
    \[
    \int_{(a,b)} \Vx{\beta}{x}(\lambda)\,\P_0(X_t \in \od x) < \infty,
    \]
    so that there is $x \in (a,b)$ such that $\Vx{\beta}{x}(\lambda)<\infty$. We conclude that for all $\lambda>0$ there is a dense set $D \subseteq [0,\infty)$ such that for all $x\in D$ one has $\Vx{\beta}{x}(\lambda)<\infty$. Then, from the scaling property,
    \[
    \Vx{\beta}{x}(\lambda) = x^{\beta+1}\Vo{\beta}(\lambda x), \quad x,\lambda>0,
    \]
    one gets $V_\beta(\lambda)<\infty$ for all $\lambda$ from the dense set in $(0,\infty)$. Since the mapping $\lambda \mapsto V_\beta(\lambda)$ is decreasing, it follows that $V_\beta(\lambda)<\infty$ for all $\lambda>0$, and the conclusion follows by one more application of the scaling property.

    $\ref{thm_moments_item_Vx} \Rightarrow \ref{thm_moments_item_ssmpp}$: since, for all $x>0$,
    \[
    \Vx{\beta}{x}(\lambda)=\int_0^\infty e^{-\lambda t}\E_x X_t^\beta \ud t<\infty,
    \]
    we immediately conclude that $\E_x X_t^\beta<\infty$ for almost all $t>0$. By arbitrary choice of $x$ we conclude $\ref{thm_moments_item_ssmpp}$.

    $\ref{thm_moments_item_ssmpp} \Rightarrow \ref{thm_moments_item_levy}$: assume that $\E_1 X_t^\beta<\infty$ for some $t>0$, or equivalently $\EE f(\xi_{\tau_t})<\infty$ with $f$ as in Lemma \ref{lem:fn_seq}. Assume first that $\beta \in (0,1]$. From \eqref{eq:moments_id} we deduce that
    \begin{align}\label{eq:26}\begin{aligned}
        \infty&>\EE \bigg[\int_0^{t \wedge \zeta} e^{-\xi_{\tau_s}}\int_{(1,\infty)} f(\xi_{\tau_s}+y)\,\nu(\od y)\ud s\bigg] \\ &\geq \EE \bigg[\int_0^{t \wedge \zeta} e^{-(1-\beta)\xi_{\tau_s}}\ind_{\{\xi_{\tau_s>0\}}} \ud s \bigg]\cdot \int_{(1,\infty)} e^{\beta y}\,\nu(\od y) \\ &= \E_1 \bigg[\int_0^{t \wedge \zeta} X_s^{\beta -1}\ind_{\{X_s>1\}}\ud s\bigg] \cdot \int_{(1,\infty)} e^{\beta y}\,\nu(\od y).\end{aligned}
    \end{align}
    Since the first factor is positive and finite, we conclude that $\int_{(1,\infty)}e^{\beta y}\,\nu(\od y)<\infty$ and \ref{thm_moments_item_levy} follows. 

    Now let $\beta \in (1,2]$. It is clear that if any (all) of the conditions \ref{thm_moments_item_ssmpp} -- \ref{thm_moments_item_V0} holds for $\beta$, it also holds for $\beta-1$, so that $\Vo{\beta-1}(\lambda)<\infty$ by the already proved part. In particular,
    \begin{equation}\label{eq:18}
    \E_1 \bigg[\int_0^{t \wedge \zeta} X_s^{\beta-1} \ud s\bigg] = \EE \bigg[ \int_0^{t \wedge \zeta} e^{(\beta-1)\xi_{\tau_s}}\ud s \bigg] \leq e^{\lambda t}\Vo{\beta-1}(\lambda)<\infty.
    \end{equation}
    Then one may observe that \eqref{eq:18} is the exact bound which is needed in \eqref{eq:14} -- \eqref{eq:16} to justify the application of the dominated convergence theorem (indeed, one then ends the estimates with $e^{-(1-\beta)\xi_{\tau_s}}\ind_{\{\xi_{\tau_s}>0\}}$ and applies \eqref{eq:18}). Thus, the equivalence \eqref{eq:moments_id} holds in the case $\beta \in (1,2]$, too. By proceeding as in \eqref{eq:26} 
    we get the statement for $\beta \in (1,2]$. Applying this procedure  $\lfloor\beta\rfloor$ times for $\beta>1$ we get the equivalence \ref{thm_moments_item_ssmpp} -- \ref{thm_moments_item_V0} in the case $\limsup_{t \to \infty}\xi_t=\infty$ $\PP$-a.s. and $\xi$ not being a compound Poisson process.
    \bigskip
    
    Now suppose that either $\lim_{t \to \infty}\xi_t=\infty$, $\lim_{t \to \infty}\xi_t=-\infty$ or is killed with independent exponential time. We will show $\ref{thm_moments_item_levy} \Rightarrow \ref{thm_moments_item_Vx} \Rightarrow \ref{thm_moments_item_ssmpp} \Rightarrow \ref{thm_moments_item_levy}$.\medskip

    $\ref{thm_moments_item_levy} \Rightarrow \ref{thm_moments_item_Vx}$: follows readily by Proposition \ref{prop:V_q_negative} and Remark \ref{remarkproposition}.

    $\ref{thm_moments_item_Vx} \Rightarrow \ref{thm_moments_item_ssmpp}$: follows as in the case $\limsup_{t \to \infty} \xi_t=\infty$ $\PP$-a.s.

    $\ref{thm_moments_item_ssmpp} \Rightarrow \ref{thm_moments_item_levy}$: follows by exactly the same recursive procedure as in the case $\limsup_{t \to \infty} \xi_t=\infty$ $\PP$-a.s, see \eqref{eq:26} and the subsequent comment.

    Thus, we have obtained the desired equivalences with \ref{thm_moments_item_ssmpp} in place of \ref{thm_moments_item_ssmp}. Since clearly $\ref{thm_moments_item_ssmp} \Rightarrow \ref{thm_moments_item_ssmpp}$, it remains to prove that any of the conditions above implies \ref{thm_moments_item_ssmp}. Indeed, fix $x>0$ and $T>0$ such that $\E_x X_T^\beta < \infty$. We claim that there is $M>0$ such that
    $\E_x X_t^\beta\leq M$ for almost all $t \in [0,T]$. Indeed, if $\beta \in (0,1]$, then one may check from the proof that the integrals $\rI_1 - \rI_4$ are uniformly bounded in $n$, $m$ and $t$ as long as $t \in [0,T]$. Furthermore, from the equivalence \eqref{eq:moments_id} we see that $\rI_5(t,m,n) \leq \lim_{n,m \to \infty} \rI_5(T,n,m)<\infty$, proving the claim. If $\beta >1$, then \eqref{eq:18} and the subsequent comment give the same assertion.
    Now, if $t$ is from the null set, then by taking a sequence $t_n \to t^+$ such that $t_n \leq T$ and
    $\E_x X_{t_n}^\beta \leq M<\infty$ for all $n \in \N$, we conclude from the right continuity and Fatou's lemma that
    \[
    \E_xX_t^\beta \leq \liminf_{n \to \infty} \E_x X_{t_n}^\beta \leq M,
    \]
    and \ref{thm_moments_item_ssmp} follows.

It remains to prove \eqref{1882}. Assume that $\beta=n\in\N$ and $\EE e^{(n+1)\xi_1}<\infty$ and fix $x>0$. By
Lemma \ref{lem:Vq} and \eqref{eq:Vq_scaling}, for all $\lambda>0$ and $k=1,\ldots,n$ one gets $\lambda \Vx{k}{x}(\lambda)=\psi(k)\Vx{k-1}{x}(\lambda)+x^{k}$. Since $\Vx{0}{x}(\lambda)=1/\lambda$, we obtain that
\begin{equation*}
\Vx{n}{x}(\lambda)=\frac{x^n}{\lambda}+
\sum_{k=1}^n\frac{\psi(n)\ldots\psi(n-k+1)}{\lambda^{k+1}}x^{n-k}.
\end{equation*}
Then \eqref{1882} is obtained for $x>0$ by inverting the Laplace transform.

Finally let us treat the case $x=0$. Let $s,t>0$ and observe that since $X_s>0$, $\P_0$-a.s.~(see for instance \eqref{eq:pssMp_ex}),  we derive from \eqref{1882} and the Markov property that,
\[
\E_0X^n_{s+t}=\E_0\big[\E_{X_s}[X_t^n]\big]=\E_0X_s^n+\sum_{k=1}^n\frac{\psi(n)\dots\psi(n-k+1)}{k!}\E_0\big[X_s^{n-k}\big]t^k.
\]
On the other hand, from the scaling property, $\E_0X^n_{s+t}=(t+s)^n\E_0X^n_{1}$ and $\E_0X^n_{s}=s^n\E_0X^n_{1}$, so that the result follows by taking $s$ to 0.
\end{proof}
\begin{remark}\label{1882killed}
    Assume $\xi$ is possibly killed with the intensity $q \geq 0$ and $\EE e^{(n+1)\xi_1}<\infty$ for some $n \in \N$. Then Lemma \ref{lem:Vq} still holds and for all $\lambda>0$ one gets $\lambda \Vx{k}{x}(\lambda)=\check{\psi}(k)\Vx{k-1}{x}(\lambda)+x^k$ for all $k =1,2,\ldots, n$, where $\check{\psi}$ is the Laplace exponent corresponding to $\check{\xi}$ which is the conservative version of $\xi$. Note that $\Vx{k}{x}(\lambda)$ is associated with possibly killed $\xi$. Then we obtain
    \[
    \Vx{n}{x}(\lambda)=\frac{x^n}{\lambda} + \sum_{k=1}^{n-1} \frac{\check{\psi}(n)\ldots\check{\psi}(n-k+1)}{\lambda^{k+1}}x^{n-k} + \check{\psi}(n)\ldots\check{\psi}(1) \frac{\Vx{0}{x}(\lambda)}{\lambda^n}.
    \]
    Thus, for all $x>0$ and $t>0$,
    \[
    \E_x X_t^n = x^n+\sum_{k=1}^{n-1} \frac{\check{\psi}(n)\ldots\check{\psi}(n-k+1)}{k!}x^{n-k}t^k+\check{\psi}(n)\ldots\check{\psi}(1)f_n(t),
    \]
    where $f$ is the inverse Laplace transform of $(0,\infty) \ni \lambda \mapsto \lambda^{-n}\Vx{0}{x}(\lambda)$. In particular, if $q=0$ and $X$ is not absorbed at $0$, i.e. $\limsup_{t \to \infty}\xi_t=\infty$, then $\Vx{0}{x}(\lambda)=1/\lambda$ and one recovers \eqref{1882}.
\end{remark}

\bibliographystyle{abbrv}
\bibliography{bib-file}
\end{document}